%% file: Ivrii_MMJ.tex
\documentclass[12pt]{article}
\usepackage[english]{babel}
\usepackage{sfmath}

\usepackage{amssymb}

\usepackage[leqno]{amsmath}

\usepackage{amsthm}

\usepackage{xr-hyper}
\usepackage[bookmarks,pdfnewwindow,plainpages=false, pdfpagelabels]{hyperref}

\hypersetup{
colorlinks=true,
linkcolor=black,
citecolor=black,
urlcolor=black,
pdfauthor={Victor Ivrii},
pdftitle={Magnetic Schr\"odinger Operator:
Geometry, Classical and Quantum Dynamics and Spectral Aymptotics},
pdfsubject={Sharp Spectral Asymptotics},
pdfkeywords={Microlocal Analysis, Magnetic Schr\"odinger Operator},
}

\usepackage{graphicx,color}
\usepackage{subfigure}

\usepackage{enumitem}
\setitemize{%
label={\textbullet},
itemindent=10pt,
leftmargin=10pt
}

\newcommand\rank{\operatorname{rank}}
\newcommand\codim{\operatorname{codim}}
\newcommand\Ker{\operatorname{Ker}}
\newcommand\Tr{\operatorname{Tr}}
\newcommand\osc{\operatorname{osc}}

\newcommand\supp{\operatorname{supp}}
\newcommand\dist{\operatorname{dist}}

\newcommand\eff{{\rm eff}}

\newcommand\bR{{\mathbb R}}
\newcommand\cE{{\mathcal E}}
\newcommand\cZ{{\mathcal Z}}
\newcommand\cV{{\mathcal V}}
\newcommand\MW{{\rm MW}}

\newcommand\const{{\rm const}}

\newcommand\Def{{\overset {\rm {def}}{\ =\ }}}

\definecolor{light}{rgb}{1,1,0.9}

\newtheorem{theorem}{Theorem}

\theoremstyle{definition}

\theoremstyle{remark}
\newtheorem{remark}[theorem]{Remark}

\begin{document}

\title{\hfill\raisebox{10pt}{\emph{\normalsize \href{http://www.math.toronto.edu/askold}{To Askold Khovanski, The Prince of Mathematics}}}\\
Magnetic Schr\"odinger Operator:
Geometry, Classical and Quantum Dynamics and Spectral Aymptotics}
\author{Victor Ivrii}
\date{\today}

\maketitle

{\abstract%
I study the Schr\"odinger operator with the strong magnetic field, considering links between geometry of magnetic field, classical and quantum dynamics  associated with operator and spectral asymptotics. In particular, I will discuss the role of short periodic trajectories.
\endabstract}

\section{Preface}
\label{sect-0}
I will consider \emph{Magnetic Schr\"odinger operator}
\begin{equation}
H= {\frac 1 2}\Bigl(\sum_{j,k}P_jg^{jk}(x)P_k -V\Bigr),\qquad P_j=hD_j-\mu A_j
\label{1}
\end{equation}
where $g^{jk}$, $A_j$, $V$ are smooth real-valued functions of $x\in \bR^d$ and
$(g^{jk})$ is a positive-definite matrix, $0<h\ll 1$ is a Planck parameter and
$\mu \gg1$ is a coupling parameter. I assume that $H$ is a self-adjoint operator.

2-dimensional magnetic Schr\"odinger is very different from 3-dimensional, all others could be close to one of these cases but are more complicated.

I am interested in the geometry of magnetic field, classical and quantum dynamics  associated with operator (\ref{1}) and spectral asymptotics
\begin{equation}
\int e(x,x,0)\psi (x)\,dx
\label{2}
\end{equation}
as $h\to +0$, $\mu\to+\infty$ where $e(x,y,\tau)$ is the Schwartz kernel of the spectral projector of $H$ and $\psi(x)$ is cut-off function. Everything is assumed to be $C^\infty$.

\section{Geometry of Magnetic Field}
\label{sect-1}
\paragraph{Magnetic Intensity}
\label{par-1-1}
Magnetic field is described by a form
\begin{equation}
\sigma = d\Bigl(\sum_{k} A_kdx_k\Bigr)={\frac 1 2}\sum_{j,k}F_{jk}dx_j\wedge d x_k
\label{3}
\end{equation}
with
\begin{equation}
F_{jk}= \partial_jV_k-\partial_kV_j.
\label{4}
\end{equation}
So $\sigma$ does not change after \emph{gauge transformation} ${\vec A}\mapsto {\vec A}+{\vec\nabla}\phi$ and this transformation does not affect other objects I am interesting in as well.

I am discussing local things and \emph{Aharonov-Bohm effect} which demonstrates that knowledge of $\sigma$, $g_{jk}$, $V$ is not sufficient to characterize spectral properties of $H$ is beyond my analysis.

\paragraph{Canonical forms} \label{par-1-2}
If $\sigma$ is of maximum rank $2r=2\lfloor d/2\rfloor$ one can reduce it locally to the \emph{Darboux canonical form}
\begin{equation}
\sigma = \sum_{1\le j\le r} dx_{2j-1}\wedge dx_{2j}.
\label{5}
\end{equation}
So, (\ref{5}) is a canonical form of $\sigma$ near \emph{generic} point for \emph{generic} $\vec A$. However situation becomes much more complicated near \emph{general} point for \emph{generic} $\vec A$. Complete results are not known. Assuming $d=2r$ and $\sigma$ is generic J.~Martinet \cite{Ma}
had shown that $\Sigma_k=\{x,\rank F(x)\le d-k\}$ are submanifolds and calculated their codimensions. In particular, $\codim \Sigma_1=1$. Moreover, $\Sigma_2=\emptyset$ as $d=2,4$ (not true for $d\ge 6$).

As $d=2$ generic 2-form $\sigma$ has a local canonical form
\begin{equation}
\sigma = x_1 dx_1\wedge dx_2, \qquad \Sigma=\{x_1=0\}.
\label{6}
\end{equation}
However, as $d=4$ not all points of $\Sigma=\Sigma_1$ are equal: $\Lambda=\{x\in \Sigma, \Ker F(x)\subset T_x\Sigma\}$ is submanifold of dimension 1. As ${\bar x}\in \Sigma\setminus \Lambda$
$\dim \bigl(\Ker F(x)\cap T_{\bar x}\Sigma\bigr)=1$ and in its vicinity
one can reduce $\sigma$ to
\begin{equation}
\sigma = x_1 dx_1\wedge dx_2 +dx_3\wedge dx_4
\label{7}
\end{equation}
while in the vicinity of ${\bar x}\in \Lambda$ canonical form is
\begin{multline}
\sigma = dx_1\wedge dx_2 - x_4dx_1\wedge dx_3 + x_3dx_1\wedge dx_4 +x_3dx_2\wedge dx_3 +x_4dx_2\wedge dx_4 +\\
2\bigl(x_1 -{\frac 1 2}(x_3^2+x_4^2)\bigr)dx_3\wedge dx_4
\label{8}
\end{multline}
(R.~Roussarie \cite{Rou}, modified by $x_2\mapsto x_2- {\frac 1 2}x_3x_4$).

\paragraph{Magnetic lines}\label{par-1-3}
\emph{Magnetic lines} are described by
\begin{equation}
{\frac {dx}{dt}}\in \Ker F(x) {\cap T_x\Sigma}
\label{9}
\end{equation}
where one can skip $ {\cap T_x\Sigma}$ without changing the substance of the definition.

As $\rank F=d$ (and thus $d$ is even) there are no magnetic lines. As $\rank F=d-1$ (and thus $d$ is odd) through each point passes exactly 1 magnetic line.

As $d=2$ and $\sigma$ is defined by (\ref{6}) the only magnetic line is  $\{x_1=0\}$.
As $d=4$ and $\sigma$ is defined by (\ref{7}) magnetic lines are straight lines $\{x_1=0, x_3=\const, x_4=\const\}$.
As $d=4$ and $\sigma$ is defined by (\ref{8}) $\Lambda=\{x_1=x_3=x_4=0\}$ and magnetic lines are helices $\{x_1=0, x_3=r\cos \theta, x_4=r\sin\theta , x_2=\const-r^2\theta /2\}$ (with $r=\const$), winging around $\Lambda$.

\paragraph{True geometry}\label{par-1-4}
From the point of view of operator $H$ simultaneous analysis of form $\sigma$ and metrics $(g^{lj})$ should be crucial, but I am not aware about any results. It appears, however, that only  eigenvalues $\pm i f_j$ and eigenspaces of matrix $(F^l_k)=(\sum_j g^{lj}F_{jk})$ are really important, and in the case of the generic magnetic field they are not very difficult to examine.

As $d=2$
\begin{equation}
f_1 = F_{12}/\sqrt g, \qquad g=\det (g^{jk})^{-1}
\label{10}
\end{equation}
while for $d=3$
\begin{equation}
f_1 = {\frac 1 2}\bigl( \sum_{j,k,l,m} g^{jk}g^{lm}F_{jl}F_{km}\bigr)^{1/2}
= \bigl( \sum_{j,k,l,m} g_{jk} F^j F^k\bigr)^{1/2}
\label{11}
\end{equation}
where $F^j={\frac 1 2}\sum _{k,l}\varepsilon^{jkl}F_{kl}$ is a vector intensity of magnetic field, $\varepsilon^{jkl}$ is an absolutely skew-symmetric tensor with $\varepsilon^{123}=1/\sqrt g$.

\section{Classical Dynamics}
\label{sect-2}
\subsection{Constant case}\label{sect-2-1}

This case has been well-known long ago. 
 \paragraph{2D case}\label{par-2-1-1}

\emph{Assume first that $g^{jk}$, $F_{jk}$ and $V$ are constant}. Then with no loss of the generality one can assume that $g^{jk}=\delta_{jk}$, skew-symmetric matrix $(F_{jk})$ is reduced to the canonical form:
\begin{equation}
F_{jk}= \left\{\begin{aligned} f_j\qquad &j=1,\dots,r, \ k=j+r\\
-f_k\qquad &j=r+1,\dots,2r, \ k=j-r\\
0\qquad &\text{otherwise}
\end{aligned}
\right.,
\label{12}
\end{equation}
$f_j>0$ and $V=0$; moreover, one can select $A_j(x)$ as linear functions.

Then as $d=2$, $f_1> 0$ classical particle described by Hamiltonian
\begin{equation}
H(x,\xi)= {\frac 1 2}\Bigl(\sum_{j,k}g^{jk}(x)\bigl(\xi_j-\mu A_j(x)\bigr) \bigl(\xi_k-\mu A_k(x)\bigr)-V(x)\Bigr)
\label{13}
\end{equation}
moves along \emph{cyclotrons} which in this case are circles of radius
$\rho_1 = ({\mu f_1})^{-1} \sqrt{2E}$ with the angular velocity $\omega_1=\mu f_1$ on energy level $\{H(x,\xi)=E\}$.

\paragraph{3D case}\label{par-2-1-2}
As $d=3$, $f_1>0$ there are a \emph{cyclotron movement} along circles of radii $\rho_1 = ({\mu f_1})^{-1} \sqrt{2E_1}$ with the angular velocity $\omega_1=\mu f_1$ and a \emph{free movement} along magnetic lines (which are straight lines along $\Ker F$) with a speed $\sqrt {2E_f}$
and energy $E$ is split into two constant arbitrary parts $E=E_1+E_f$.

 \paragraph{Multidimensional case}\label{par-2-1-3}

Multidimensional case with $d=2r=\rank F$ is a combination of 2D cases:
there are $r$ cyclotron movements with angular velocities $\omega_k=\mu f_k$
and radii $\rho_k = ({\mu f_k})^{-1} \sqrt{2E_k}$
where energy $E$ is split into $r$ constant arbitrary parts $E=E_1+E_2+\dots +E_r$. The exact nature of the trajectories depends on the comeasurability of $f_1,\dots, f_r$.

As $d>2r=\rank F$ in addition to the cyclotronic movements described above appears a free movement along any constant direction ${\vec v}\in \Ker F$ with a speed $\sqrt{2E_f}$ where energy $E$ is split into $r+1$ constant arbitrary parts $E=E_1+E_2+\dots +E_r+E_f$.

This difference between cases $d=2r=\rank F$ and $d>2r=\rank F$ will be traced through the whole paper.

\subsection{Full rank case}
\label{sect-2-2}

\emph{Assume now only that $d=\rank F$} (see f.e \cite{IRO3,IRO4}. In addition,  \emph{assume temporarily that $F_{jk}$ and $g^{jk}$ are constant but potential $V(x)$ is linear}.
Then cyclotronic movement(s) is combined with the \emph{magnetic drift} described by equation
\begin{equation}
{\frac {d x_j}{dt}}= (2\mu)^{-1} \sum_k \Phi^{jk}\partial_k V
\label{14}
\end{equation}
where $(\Phi^{jk})=(F_{jk})^{-1}$.

As $d=2$ it will be movement along cycloid
and multidimensional movement will be combination of those.

Not assuming anymore that $V$ is linear we get a bit more complicated picture:

\begin{itemize}

\item Equation (\ref{14}) holds modulo $O(\mu^{-2})$; modulo error $O(\mu^{-2}t)$;

\item As $d=2$ cycloid is replaced by a more complicated curve drifting along $V=\const$ and thus cyclotron radius $\rho=(\mu f_1)^{-1}\sqrt{2E+V}$ would be
preserved;

\item In higher dimensions all cyclotron radii are preserved as well.

\end{itemize}

Without assumption that $g^{jk}$ annd $F_{jk}$ are constant picture becomes even more complicated:

\begin{itemize}

\item As $d=2$ cycloid is replaced by a more complicated curve drifting along $f^{-1}(V+2E)=\const$ (thus preserving \emph{angular momentum}
$\omega_1\rho_1^2$ according to equation
\begin{equation}
{\frac {d x}{dt}}= (2\mu )^{-1} \bigl(\nabla f^{-1}(V+2E)\bigr)^\perp
\label{15}
\end{equation}
where $^\perp$ means clockwise rotation by $\pi/2$ assuming that at point in question $g^{jk}=\delta_{jk}$;

\item In higher dimensions (at least as \emph{non-resonance conditions $f_j\ne f_k$ $\forall j\ne k$} and \emph{$f_j\ne f_k+f_l$ $\forall j,k,l$} are fulfilled) one can split potential $V=V_1+\dots +V_k$ so that similar equations hold in each eigenspace of $(F^j_k)$ and both separate energies and angular momenta are (almost) preserved.
\end{itemize}

\subsection{3D case}
\label{sect-2-3}
As $d>2r=\rank F$ the free movement is the main source of the spatial displacement and the most interesting case is $2r=d-1$ and especially $d=3$, $r=1$.

In this case  the \emph{magnetic angular momentum} $M$ is (almost) preserved; thus \emph{kinetic energy} of magnetic rotation is ${\frac 1 2}f^{-1}M^2$; therefore in the coordinate system such that $g^{1j}=\delta_{1j}$ the free movement is described by 1D Hamiltonian
\begin{equation}
H_1 (x_1,\xi_1; x',M)= {\frac 1 2}\xi_1^2 -{\frac 1 2} V_\eff
\label{16}
\end{equation}
with effective potential $V_\eff (x_1,x')=V - f^{-1}M^2$, $x=(x_1,x')$.

Thus the particle does not necessarily run the whole magnetic line and the helix winging around it does not necessarily have constant the step or radius.
Effect of the magnetic drift is rather minor.

\subsection{2D case: variable rank}
\label{sect-2-4}
Situation becomes really complicated for variable $\rank F$. I am going to consider only $d=2,4$ and a generic magnetic form $\sigma$. Let me start from the model Hamiltonian as $d=2$:
\begin{equation}
H^0={\frac 1 2}\Bigl(\xi_1^2 + (\xi_2-\mu x_1^\nu /\nu )^2-1\Bigr);
\label{17}
\end{equation}
the drift equation is
\begin{equation}
{\frac {dx_1}{dt}}=0, \quad {\frac {dx_2}{dt}}=
{\frac 1 2}(\nu-1)\mu ^{-1}x_1^{-\nu}
\label{18}
\end{equation}
and for $|x_1|\gg {\bar\gamma}=\mu^{-1/(\nu+1)}$ gives a proper description of the picture.

For the model Hamiltonian (\ref{17}) with $\mu=1$ (otherwise one can scale $x_1\mapsto \mu^{1/2}x_1$, $\xi_2\mapsto \mu k$) we can consider also 1-dimensional movement along $x_1$ with potential
\begin{equation}
\cV(x_1;k) =1-(k-x_1^\nu/\nu)^2, \qquad k=\xi_2;
\label{19}
\end{equation}
Then for odd $\nu$
\begin{itemize}
\item $\cV$ is one-well potential;
\item As $k=\pm 1$ one of its extremes is $0$ where ${\frac {d\cV}{dx_1}}(0)=0$;
\item Well is more to the right/left from 0 as $\pm k>0$; as $k=0$ well becomes symmetric.
\end{itemize}
%%%%%%%%%%%%%%%%%%%%%%%%%%%%%%%%%%%%%%%%%%%%%%%%%%%
% Figure fig-1. Note that I use package subfigure
%%%%%%%%%%%%%%%%%%%%%%%%%%%%%%%%%%%%%%%%%%%%%%%%%%%
\begin{figure}[h!]
\centering
\subfigure[$k> 1$;]{
\includegraphics[width=.3\textwidth]{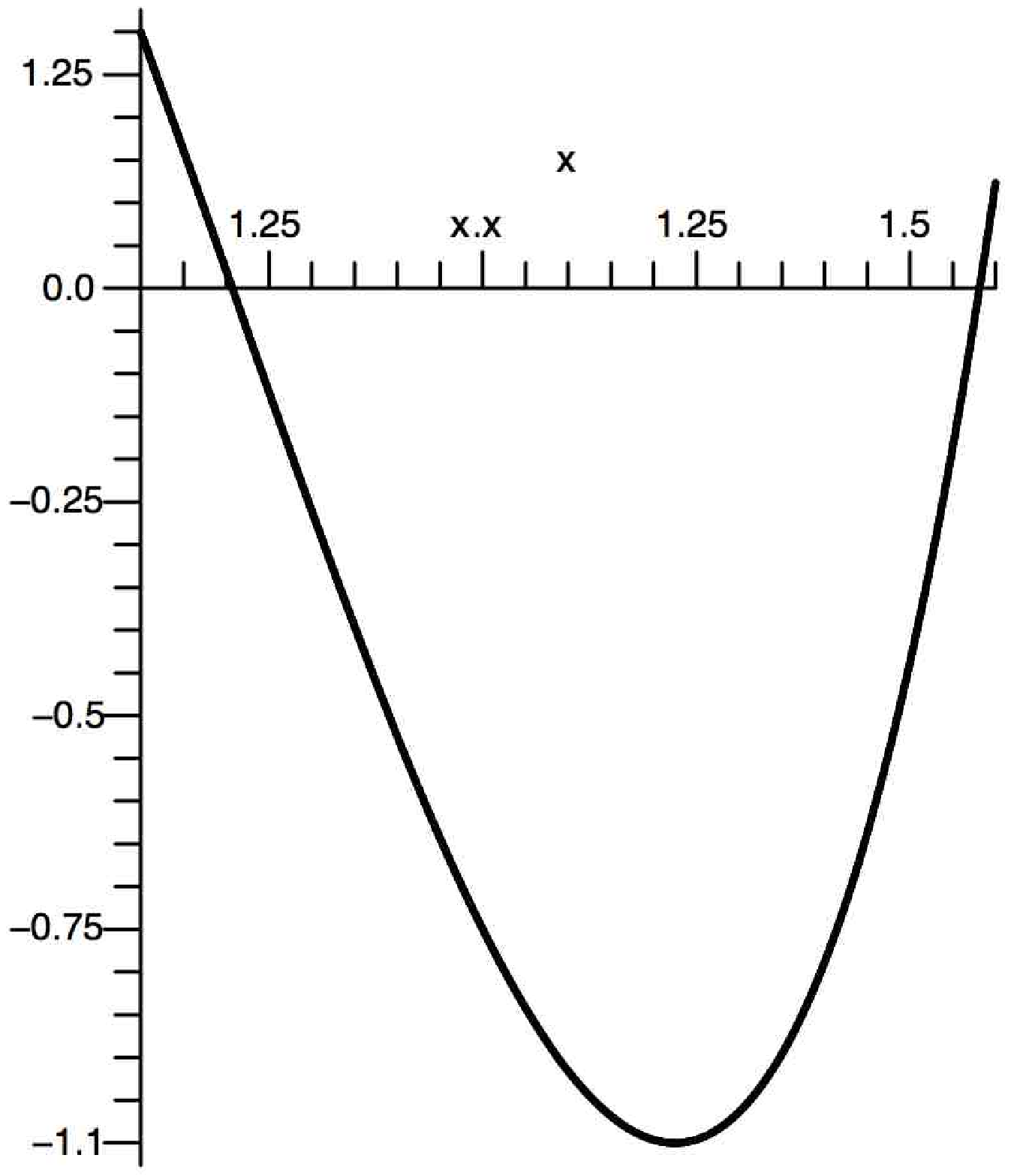}}
\subfigure[ $0<k< 1$;]{
\includegraphics[width=.3\textwidth]{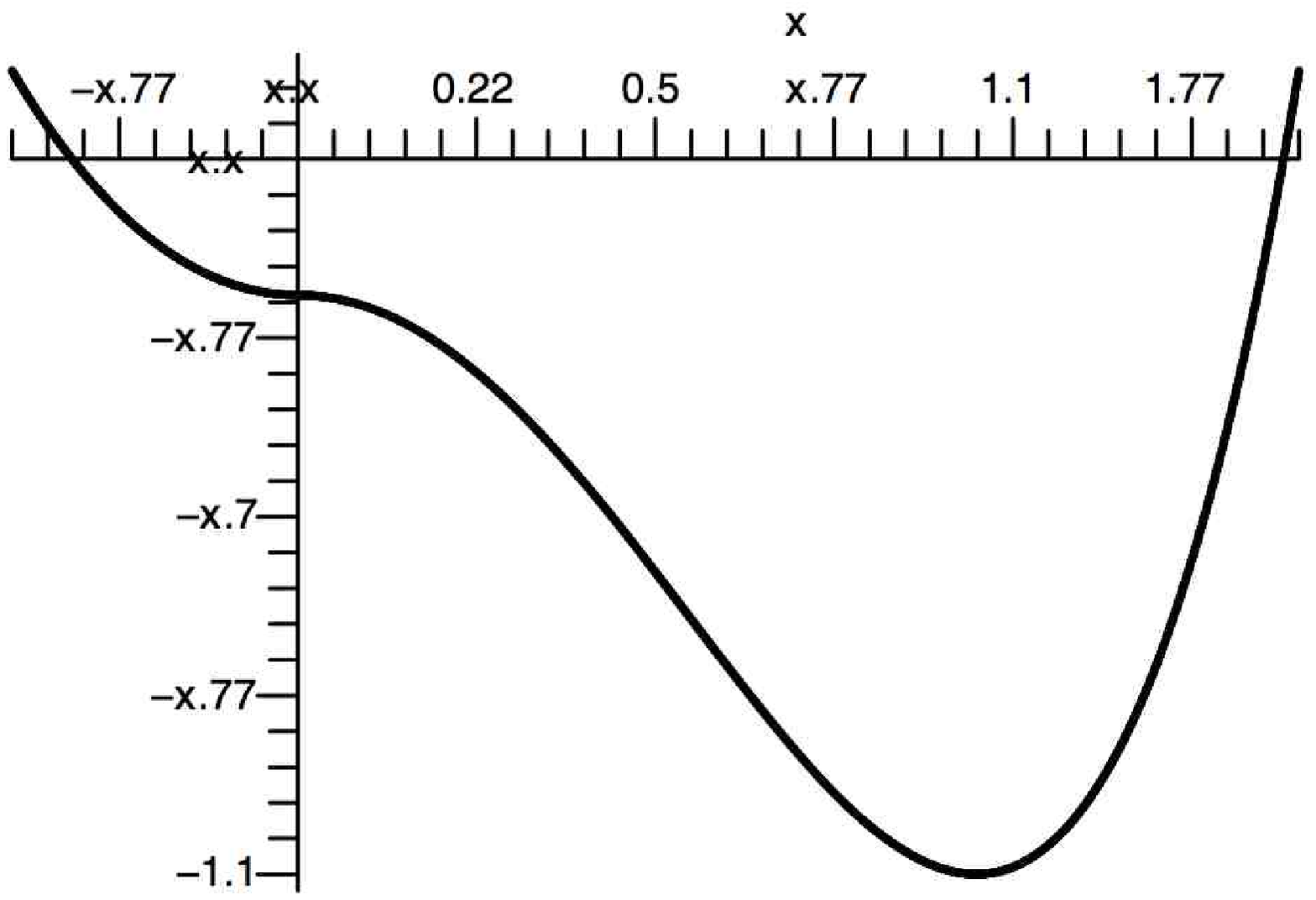}}
\subfigure[$k= 1$;]{
\includegraphics[width=.3\textwidth]{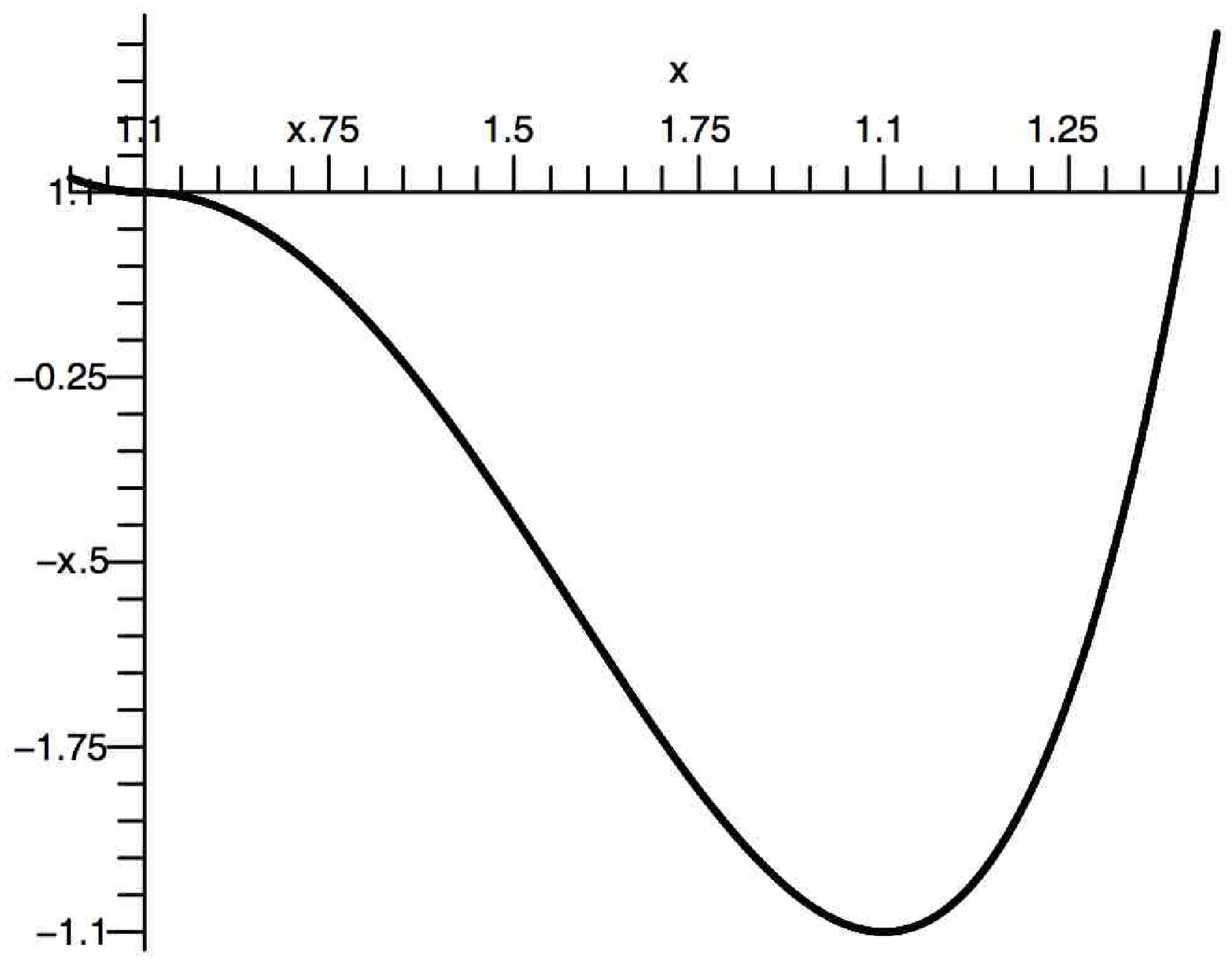}}
\caption{Effective potential for odd $\nu$.}
\label{fig-1}
\end{figure}
On the other hand, for even $\nu$ potential is always symmetric and
\begin{itemize}
\item We have two-well potential with the central bump above surface if $k >1$ \item and below it as $0<k<1$:
\item As $k=\pm 1$ one of its extremes is $0$ where ${\frac {d\cV}{dx_1}}(0)=0$;
\item Well is more to the right/left from 0 as $\pm k>0$; as $k=0$ well becomes symmetric.
\end{itemize}
%%%%%%%%%%%%%%%%%%%%%%%%%%%%%%%%%%%%%%%%%%%%%%%%%%%
% Figure fig-2. Note that I use package subfigure
%%%%%%%%%%%%%%%%%%%%%%%%%%%%%%%%%%%%%%%%%%%%%%%%%%%
\vskip-10pt
\begin{figure}[h!]
\centering
\subfigure[$k> 1$;]{
\includegraphics[width=.23\textwidth]{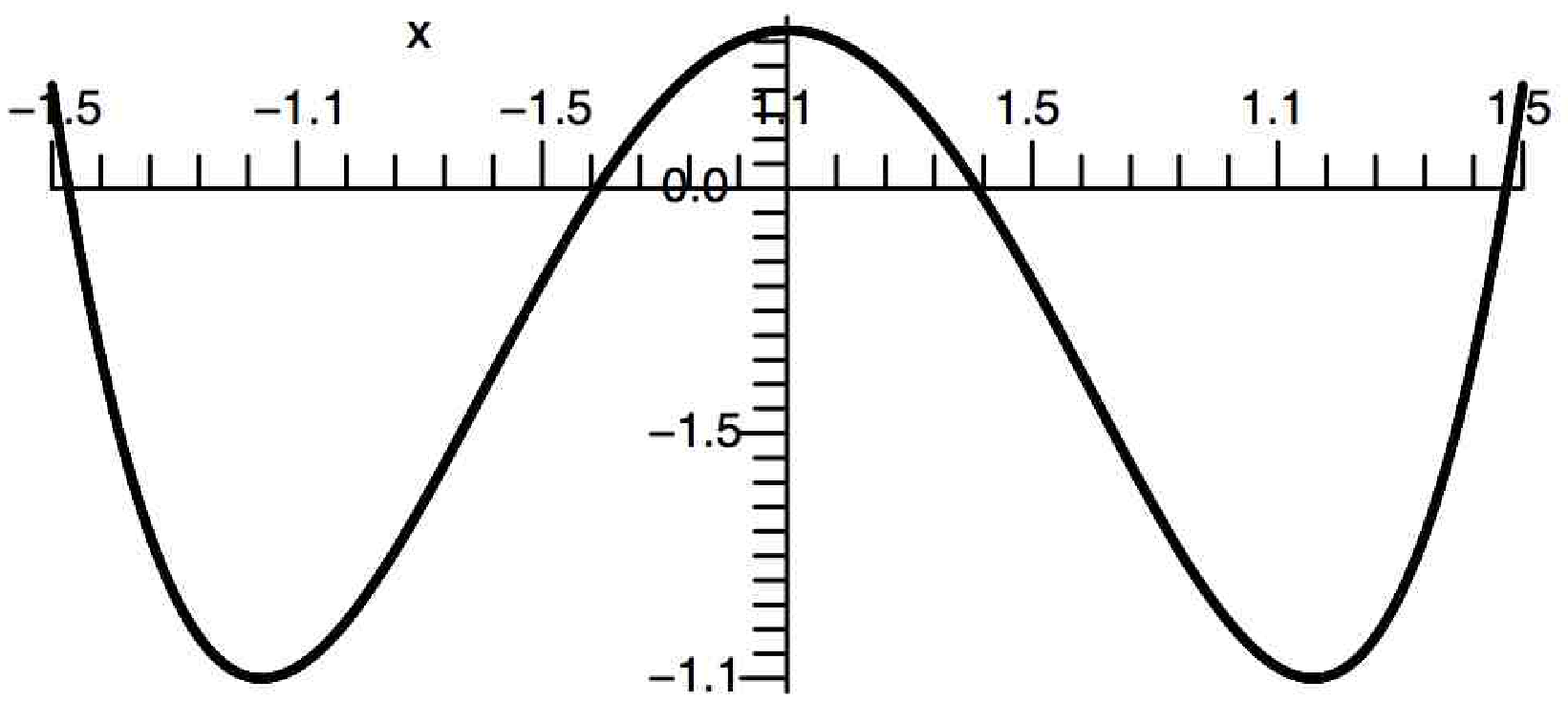}}
\subfigure[$k=1$;]{
\includegraphics[width=.23\textwidth]{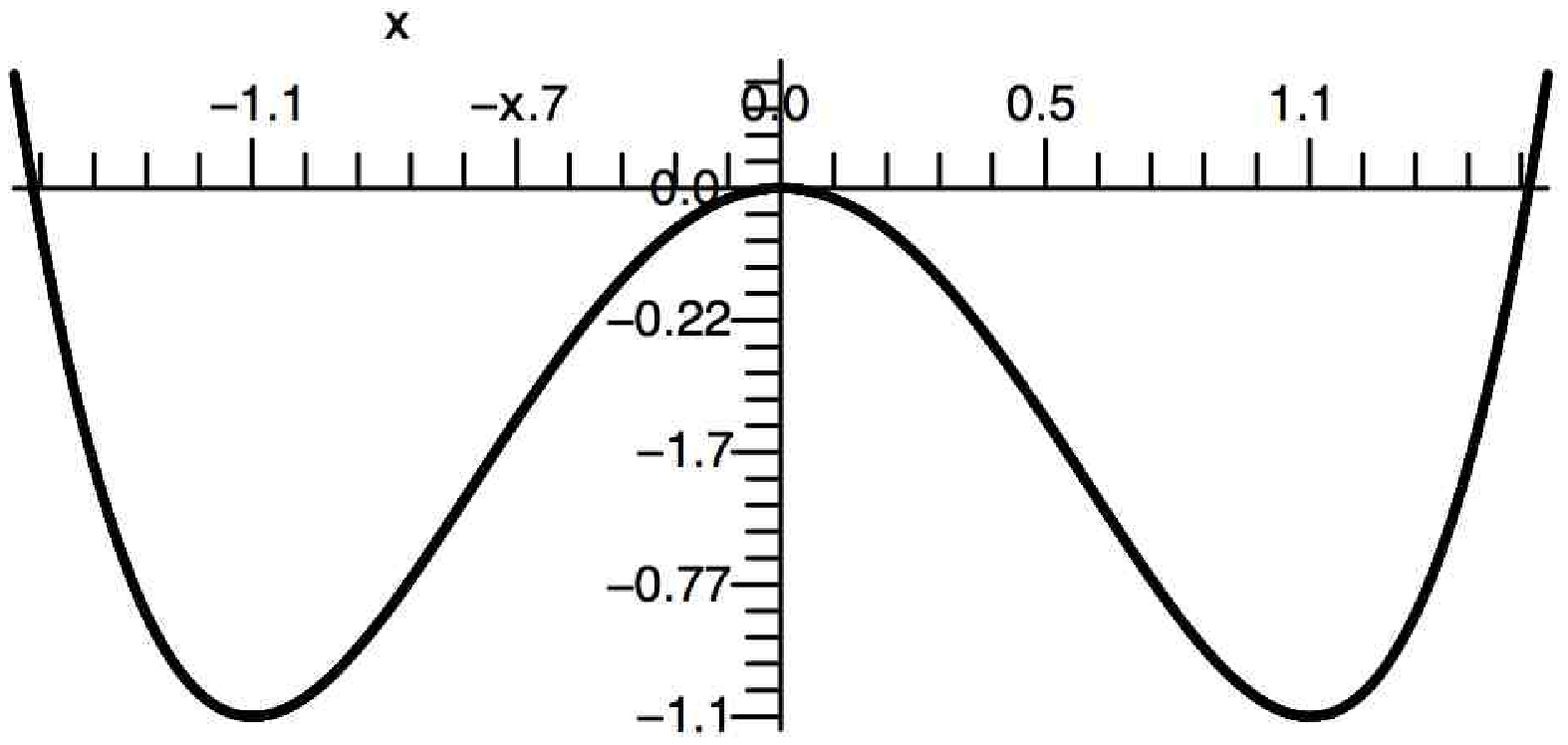}}
\subfigure[$0<k< 1$;]{
\includegraphics[width=.23\textwidth]{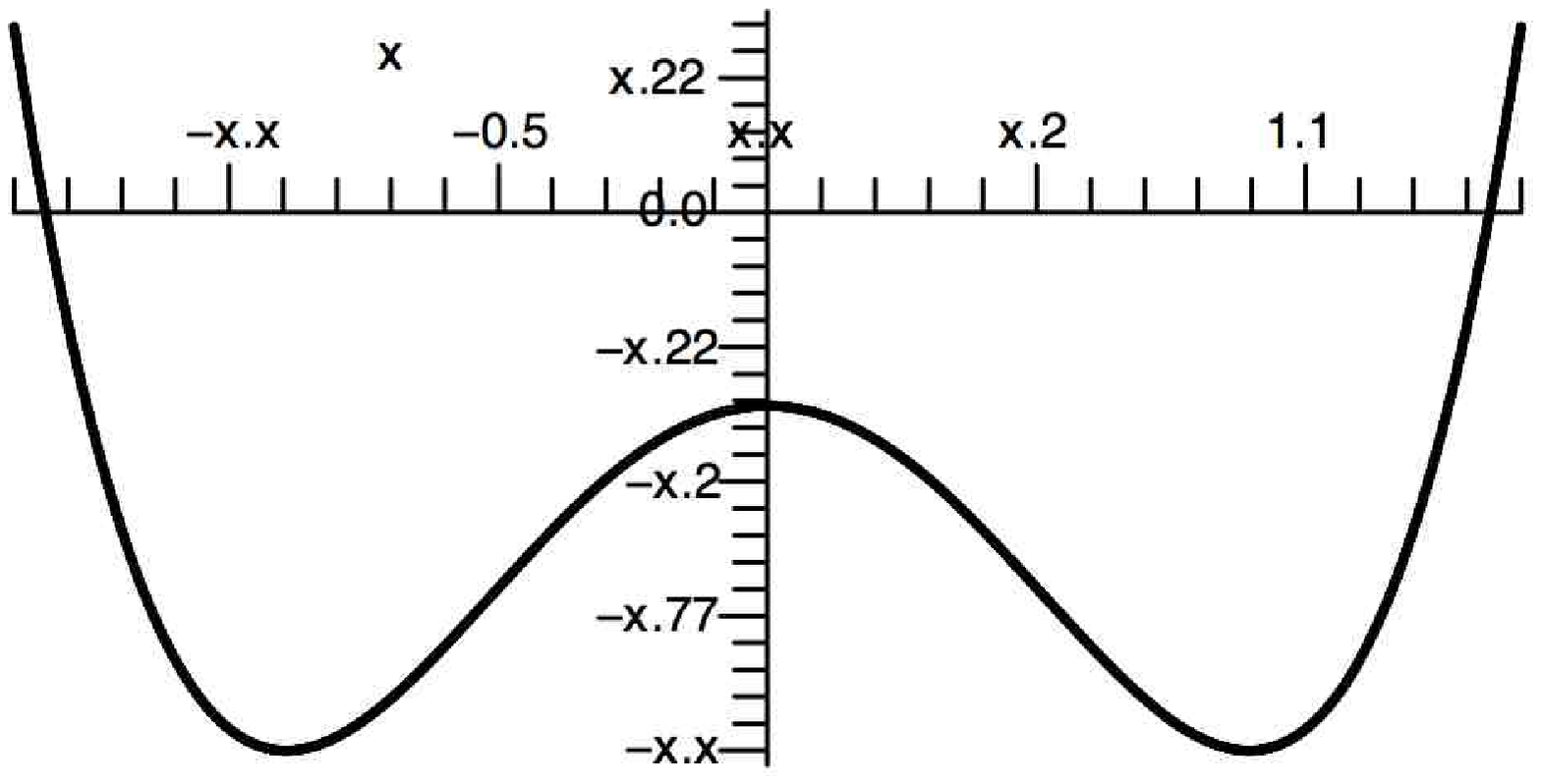}}
\subfigure[$-1<k\le 0$;]{
\includegraphics[width=.23\textwidth]{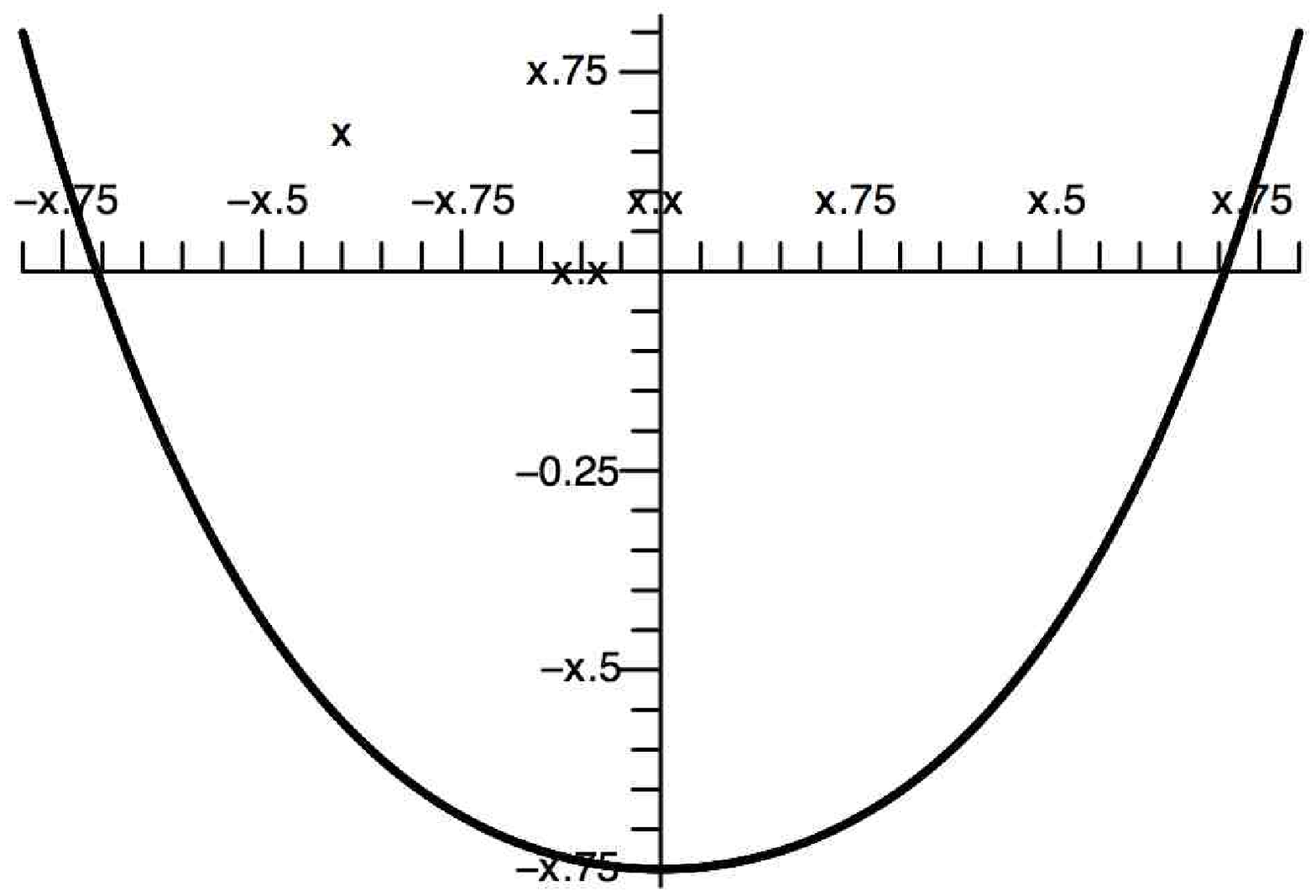}}
\caption{Effective potential for even $\nu$.}
\label{fig-2}
\end{figure}

Let us consider trajectories \emph{on the energy level $0$}.
From the analysis of the potential one can  conclude that for $k\ne \pm 1$ the movement along $x_1$ is periodic with the period
\begin{equation}
T(k)=2\int_{x^-_1(k)}^{x^+_1(k)} {\frac {dx_1}{\sqrt {2\cV(x_1;k)}}};
\label{20}
\end{equation}
however one needs to analyze the increment of $x_2$ during this period
\begin{equation}
I(k)=2\int_{x^-_1(k)}^{x^+_1(k)} {\frac {(k-x_1^\nu/\nu)dx_1}{\sqrt {2\cV(x_1;k)}}}.
\label{21}\end{equation}
One can prove that $I(k)\gtrless 0$ as $k\gtrless k^*$ with $0<k^*<1$ for even $\nu$ and $k^*=0$ for odd $\nu$. In particular, $k^*\approx 0.65$ for $\nu=2$.
Further, $I(k)\asymp (k-k^*)$ as $k \approx k^*$.

\begin{figure}[h!]
\centering
\subfigure[$k\gg 1$; as $x_1>0$ trajectory moves up and rotates clockwise]{
\includegraphics[width=0.4\textwidth]{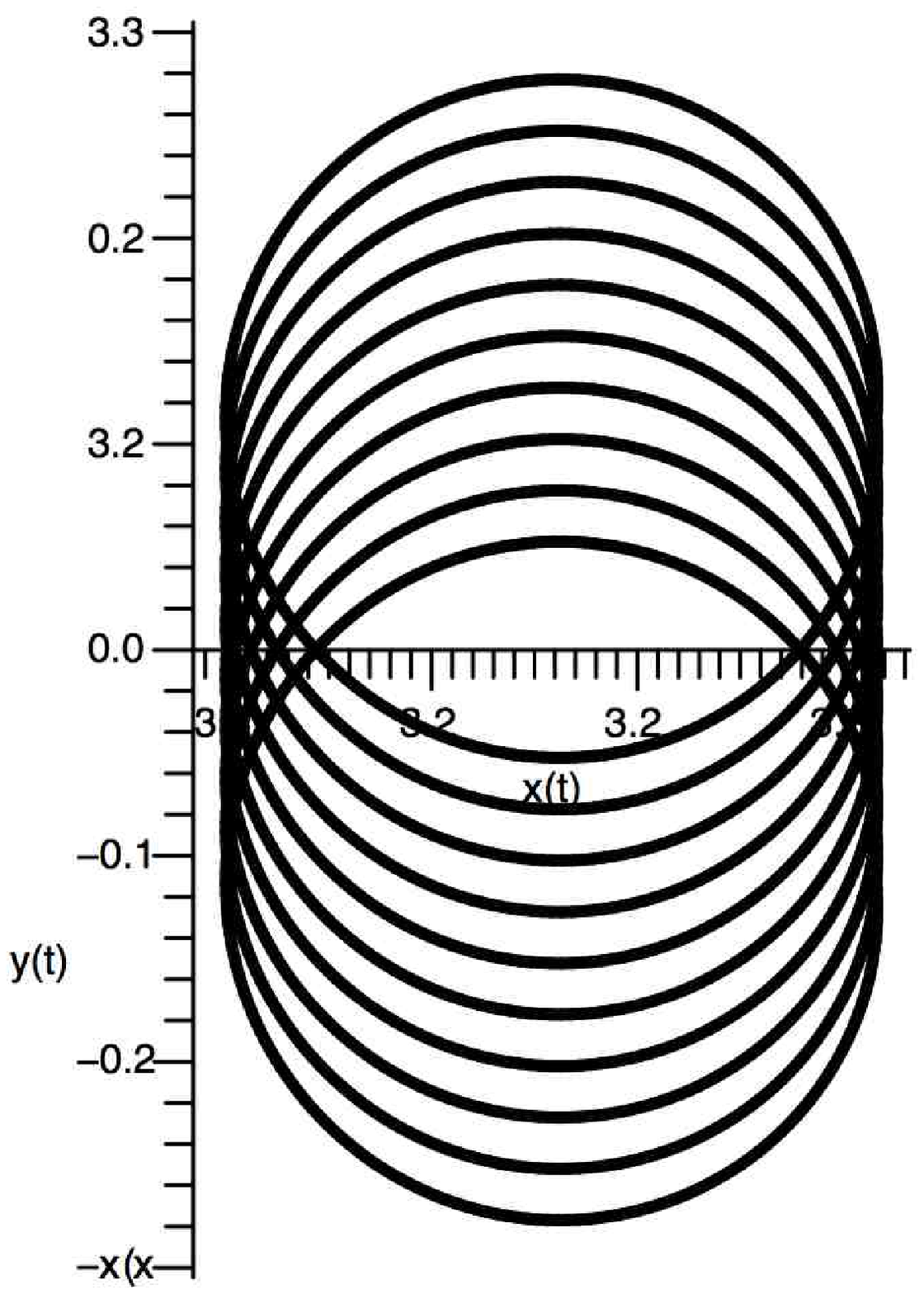}}
\subfigure[$k$ decreases, still $k>1$. Trajectory becomes less tight; actual size of cyclotrons increases;]{
\includegraphics[width=0.4\textwidth]{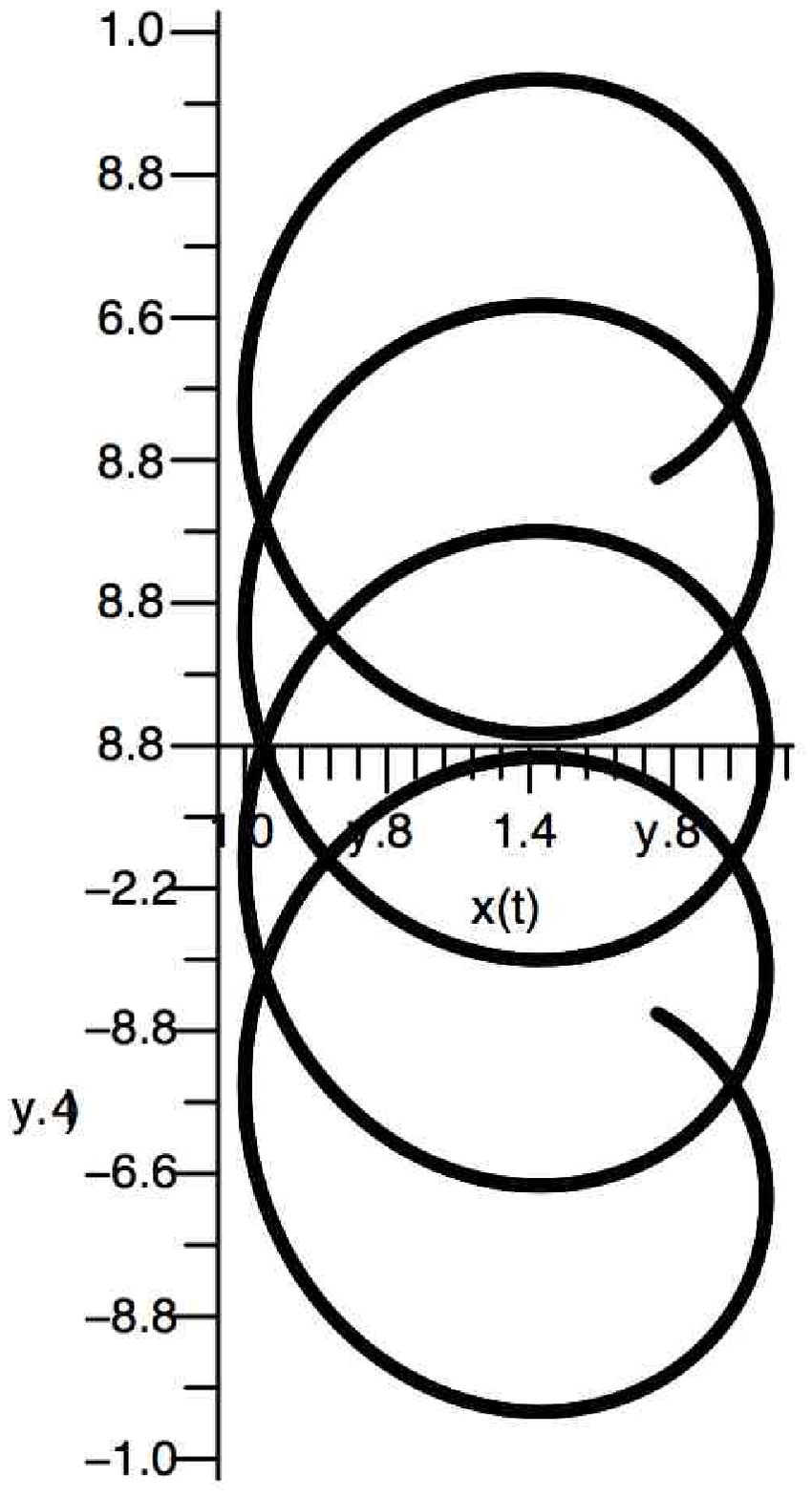}}
\subfigure[$k$ further decreases, still $k>1$;]{
\includegraphics[width=0.4\textwidth]{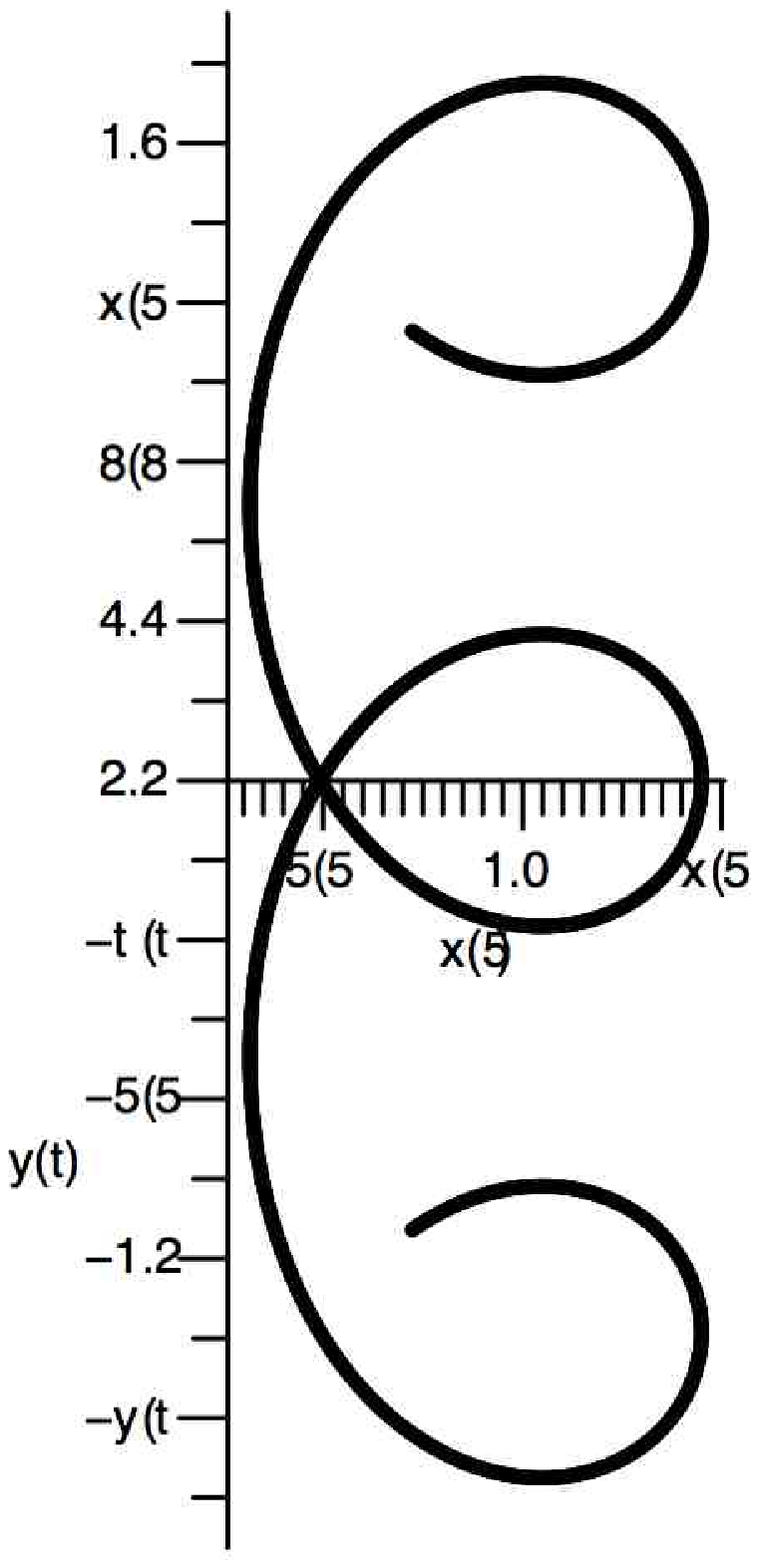}}
\subfigure[$k=1$. Trajectory contains just one cyclotron.]{
\includegraphics[width=0.4\textwidth]{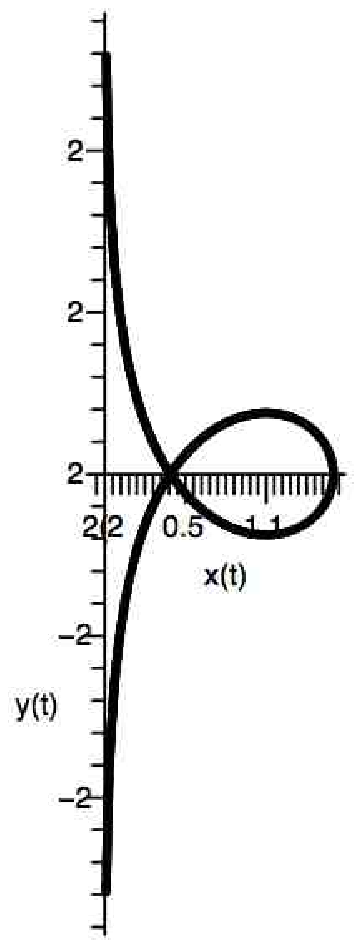}}
\caption{Movements in the outer zone}
\label{fig-3}
\end{figure}

\begin{figure}[p]
\centering
\subfigure[$k<1$ slightly]{
\includegraphics[width=0.4\textwidth]{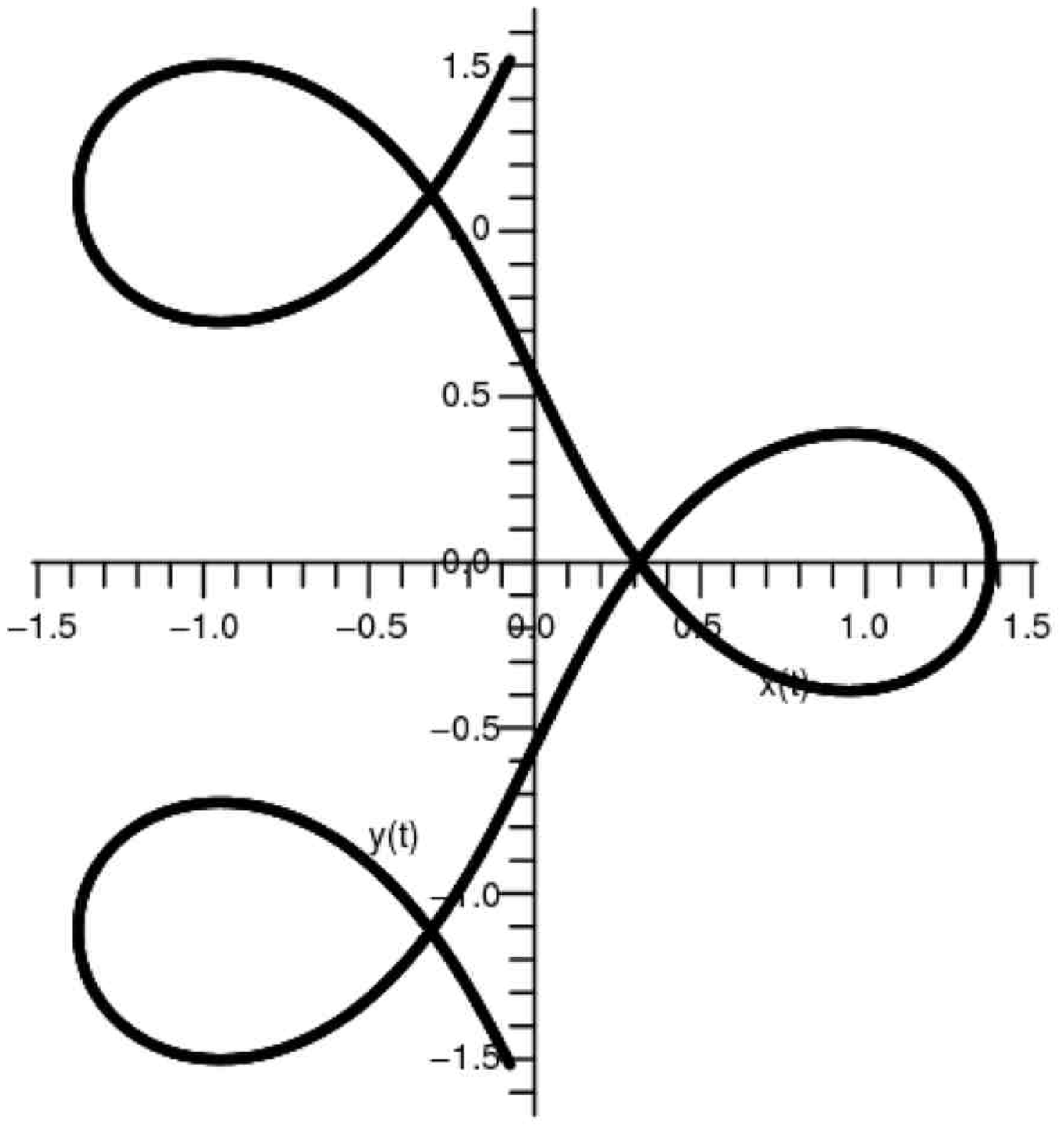}}
\subfigure[$k$ further decays but still larger than $k^*$. Drift slows down]{
\includegraphics[width=0.4\textwidth]{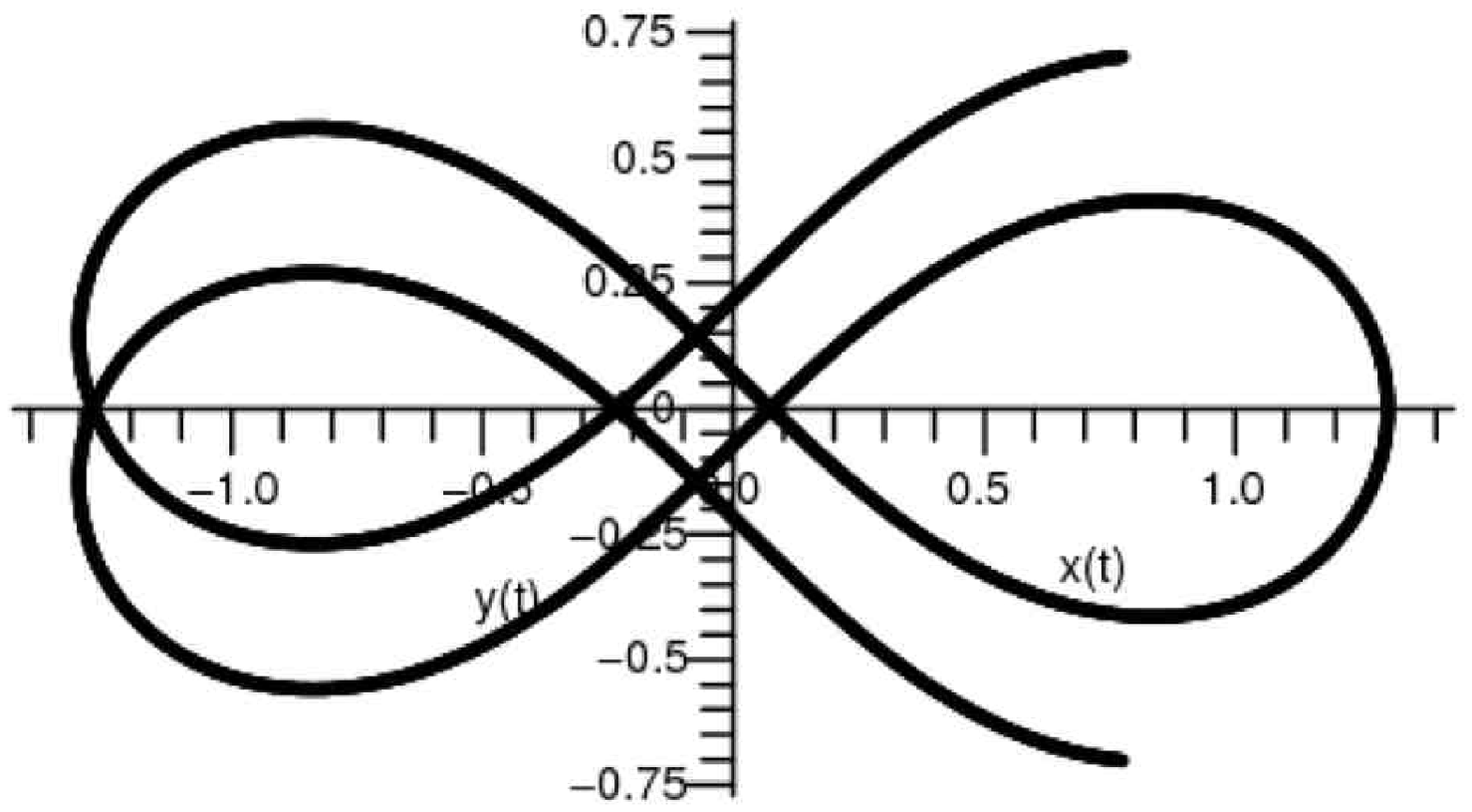}}
\subfigure[$k=k^*$. No drift; trajectory becomes periodic]{
\includegraphics[width=0.4\textwidth]{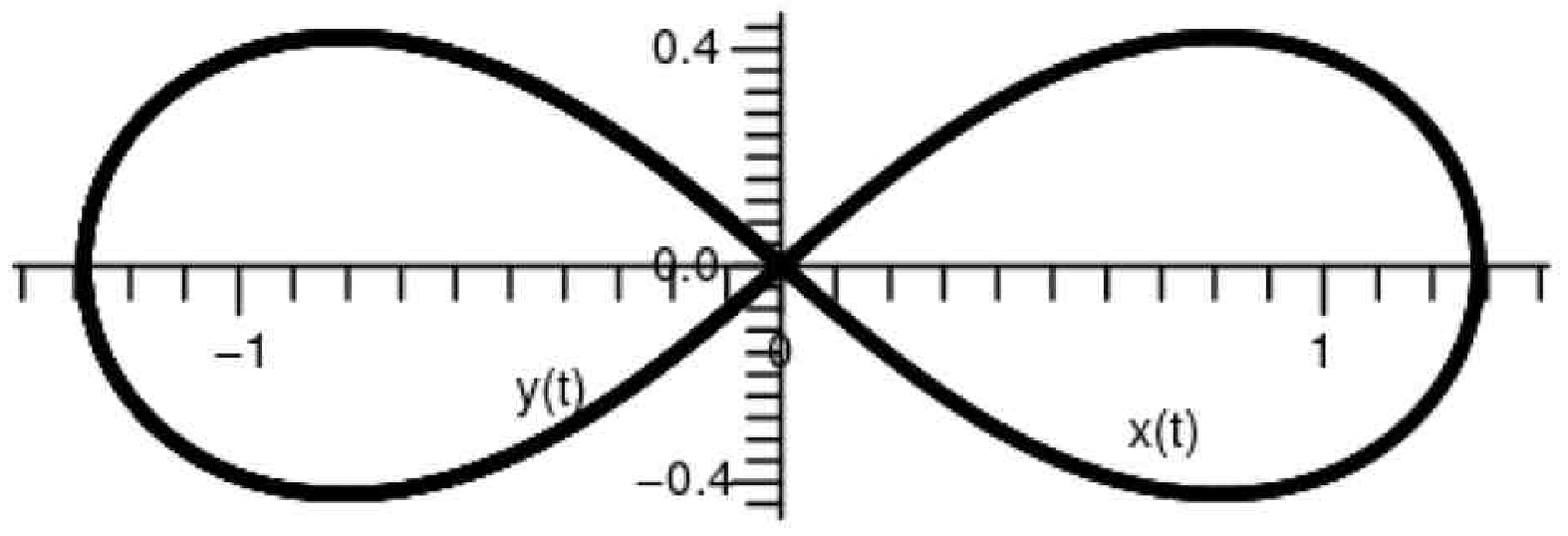}}
\subfigure[$k< k^*$. Drift now is down!]{
\includegraphics[width=0.4\textwidth]{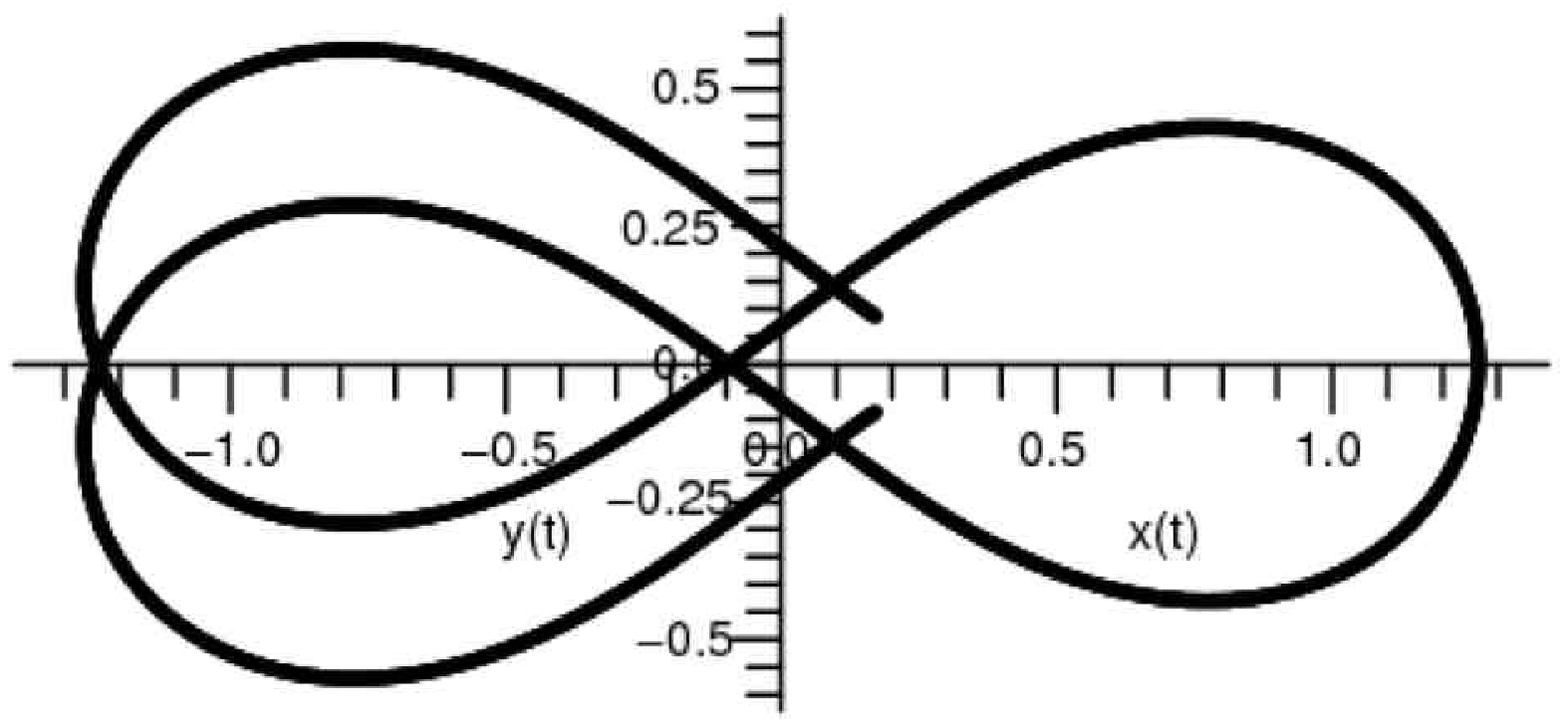}}
\subfigure[$k$ decays further. Drift down accelerates.]{
\includegraphics[width=0.4\textwidth]{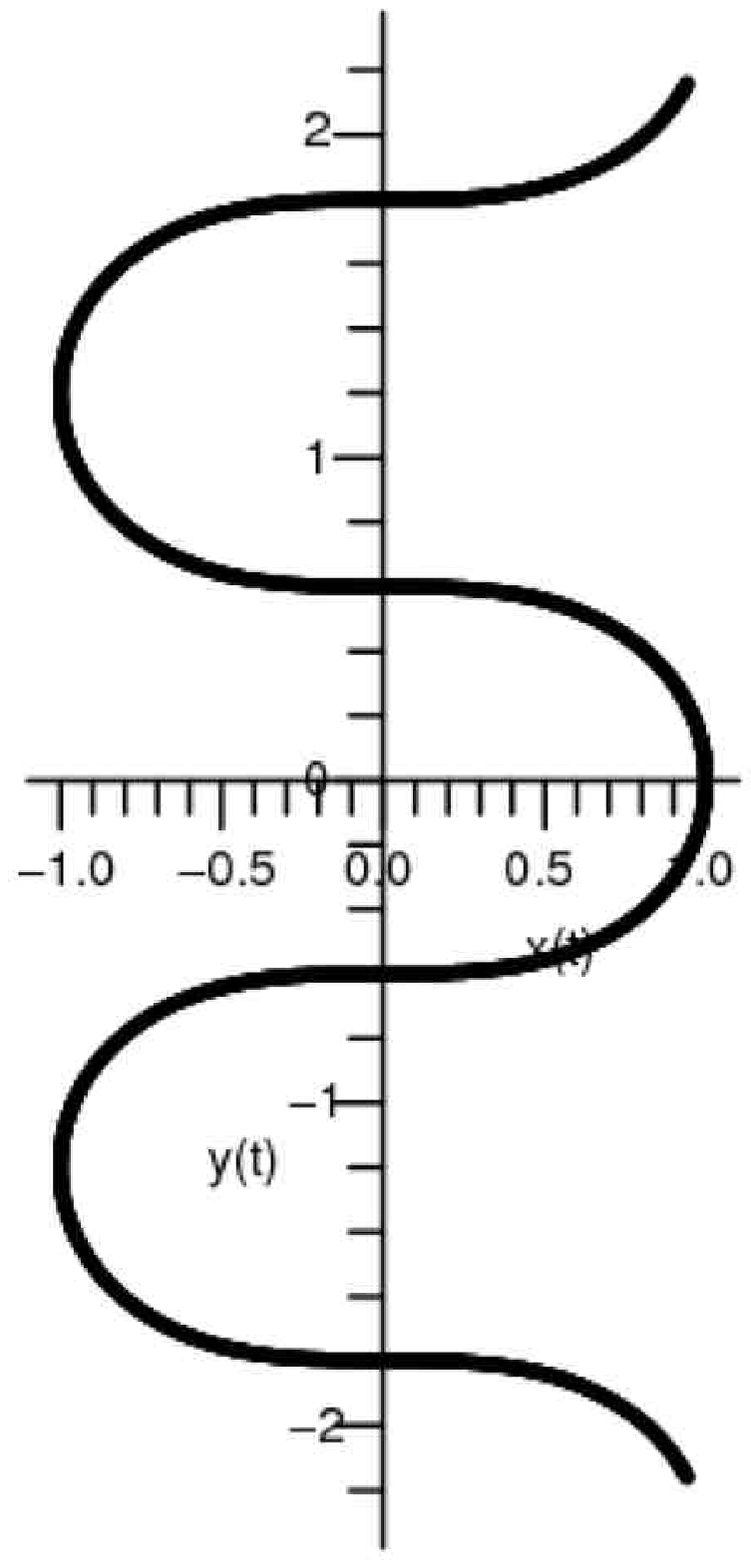}}
\subfigure[and further; as $k=-1$ we have just straight line down]{
\includegraphics[width=0.4\textwidth]{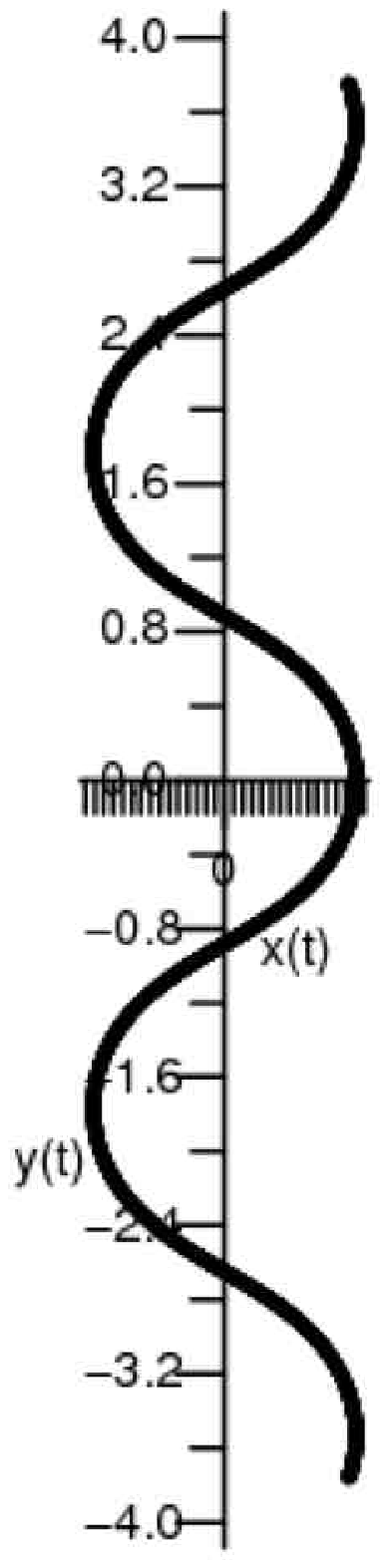}}
\caption{Even $\nu$; movements in the inner zone}
\label{fig-4}
\end{figure}

\begin{figure}[h!]
\centering
\subfigure[$k<1$ slightly. Drift is up and the fastest]{
\includegraphics[width=0.4\textwidth]{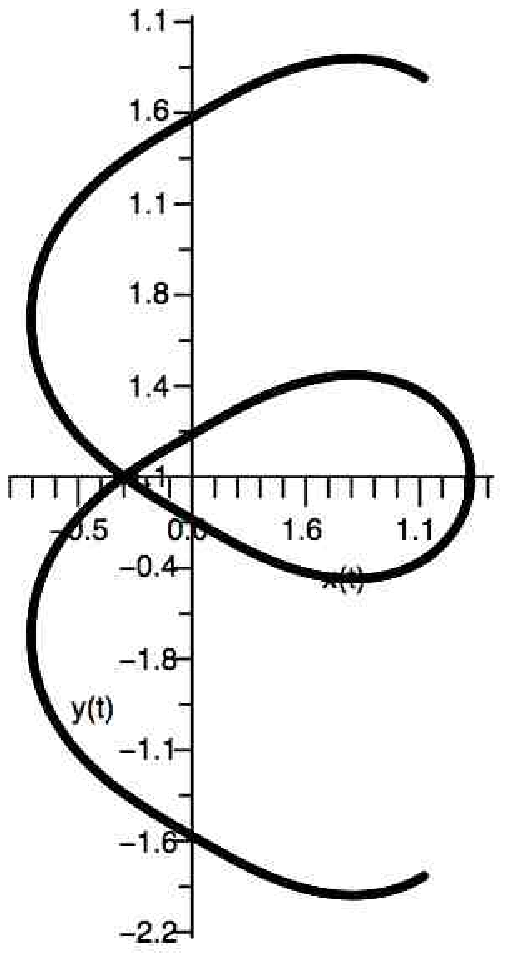}}
\subfigure[$k$ decays further. Drift up slows down further]{
\includegraphics[width=0.25\textwidth]{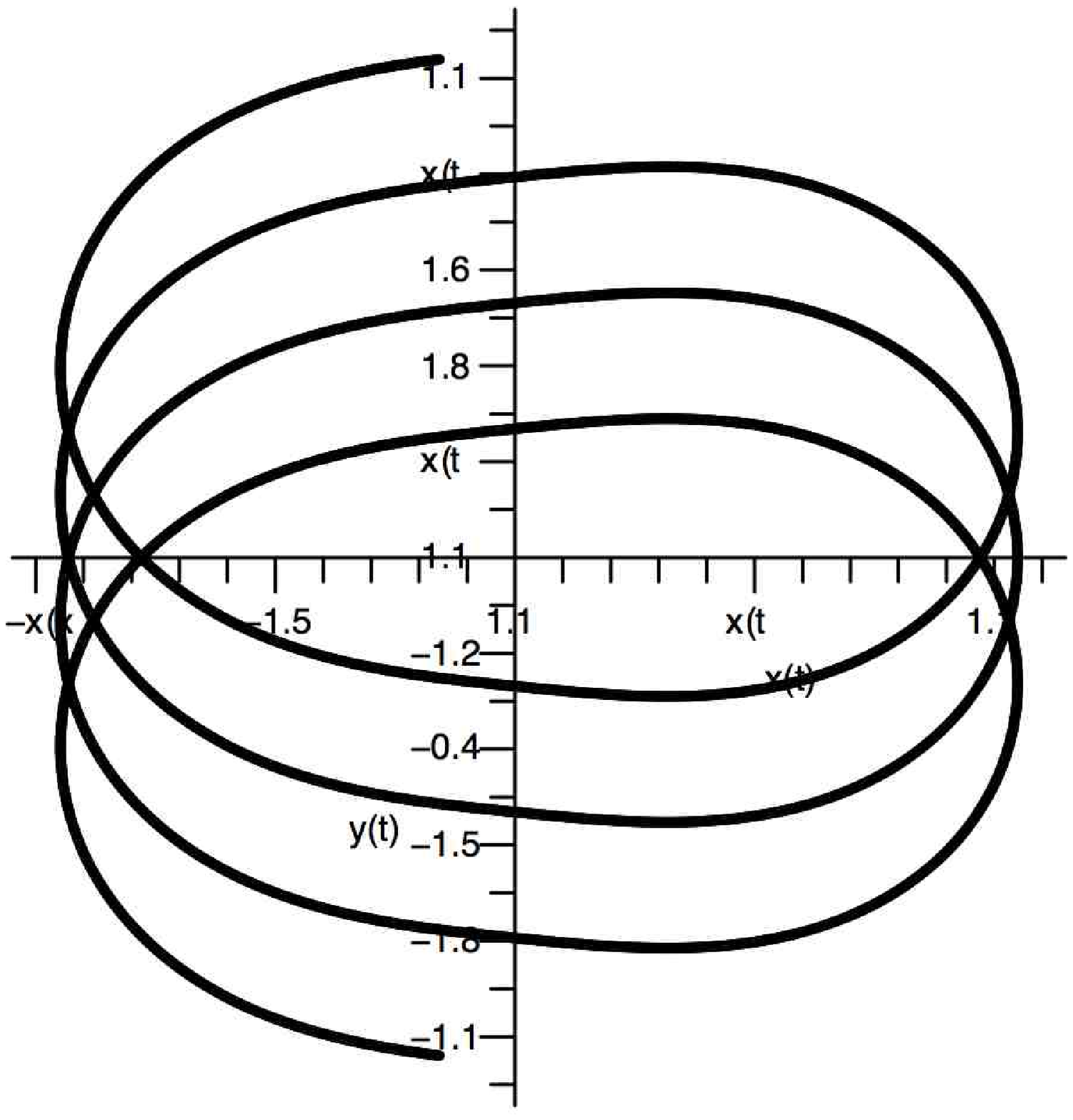}}
\subfigure[$k=0$. No drift. Trajectory is periodic]{
\includegraphics[width=0.25\textwidth]{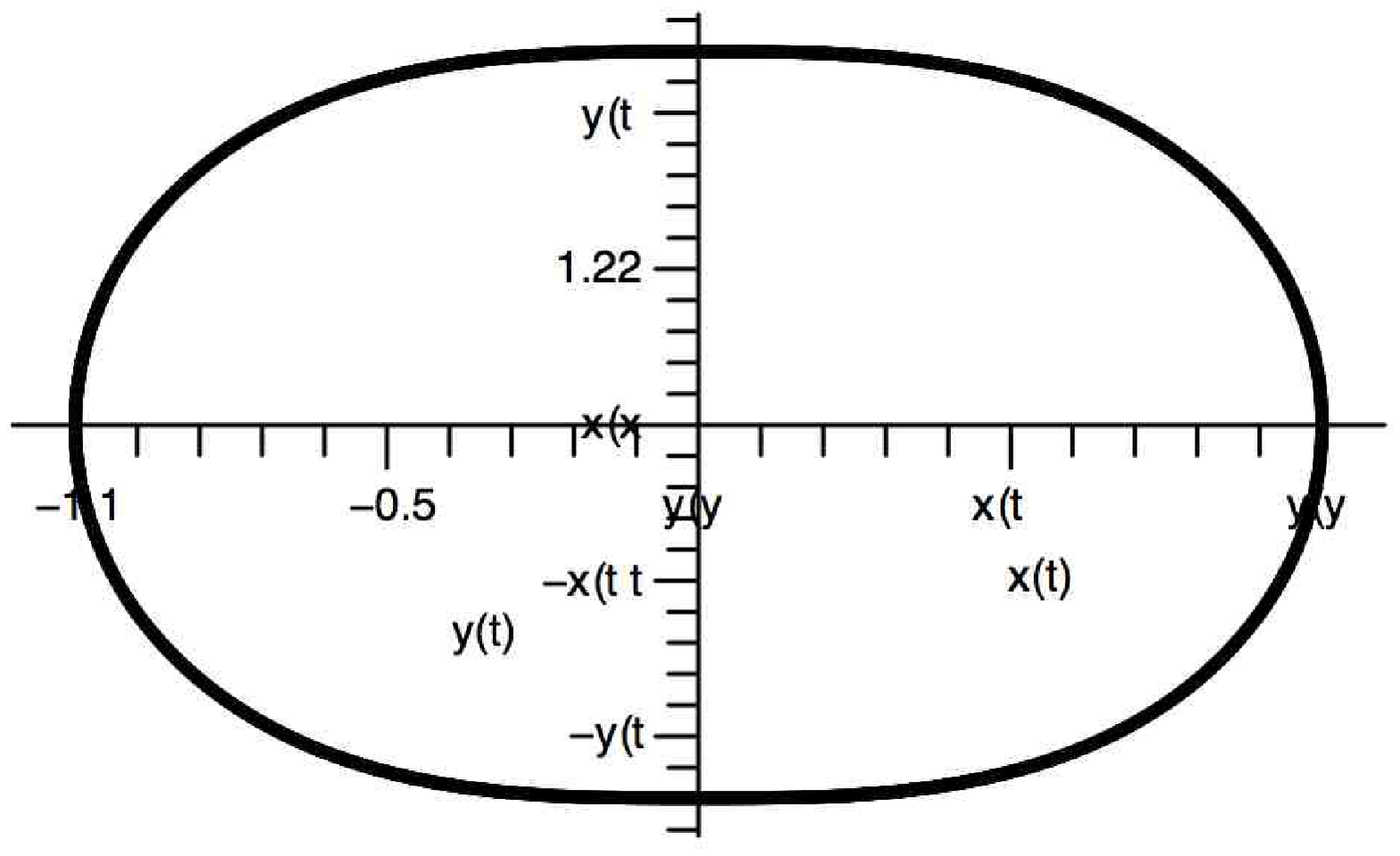}}
\caption{Odd $\nu$; movements in the inner zone}
\label{fig-5}
\end{figure}

\begin{figure}[h!]
\centering
\subfigure[Even $\nu$]{
\includegraphics[width=0.35\textwidth]{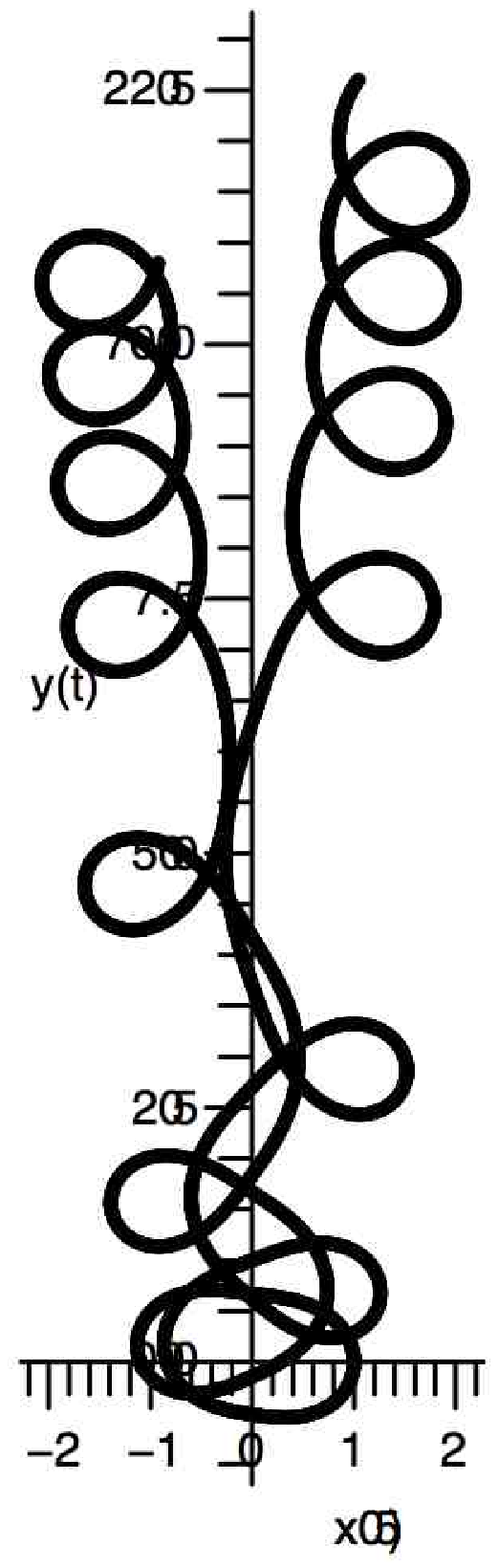}}
\subfigure[Odd $\nu$]{
\includegraphics[width=0.35\textwidth]{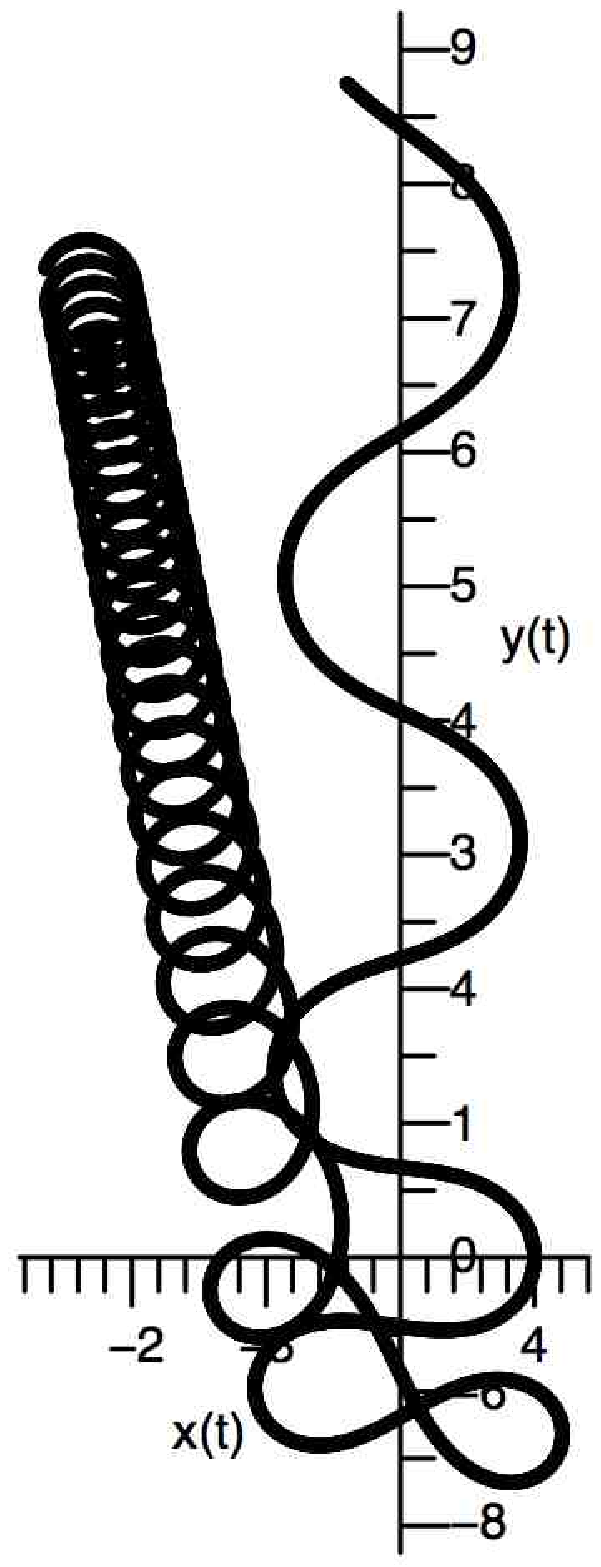}}
\caption{Breaking periodic trajectories by a linear potential}
\label{fig-6}
\end{figure}

On figures \ref{fig-3}--\ref{fig-5} are shown trajectories on $(x_1,x_2)$-plane plotted by Maple in the outer zone (these trajectories have mirror-symmetric or central-symmetric for even or odd $\nu$ in zone $x_1<0$), in the inner zone for even $\nu$ and in the inner zone for odd $\nu$ respectively. For the spectral asymptotics periodic trajectories are very important, especially the short ones. Periodic trajectories shown above are \emph{very unstable} and taking $V=1-\alpha x_1$ instead of $x_1$ breaks them down
(figure \ref{fig-6}).
%%%%%%%%%%%%%%%%%%%%%%%%%%%%%%%%%%%%%%%%%%%%%%%%%%%%%%%%%%%%%%%%%%%%%%
% Figures fig-3, fig-4, fig-5
% Clearly the best place for figures would be here but since breaking them 
% is less than desirable, compromise ``place them on the next pages'' seems
% to be the best. Note that I use package subfigure
%%%%%%%%%%%%%%%%%%%%%%%%%%%%%%%%%%%%%%%%%%%%%%%%%%%%%%%%%%%%%%%%%%%%%%%
%%%%%%%%%%%%%%%%%%%
%%%%%%%%%%%%%%%%%%%%%%%%%%
%%%%%%%%%%%%%%%%%%%%%%%%%%%%%%%%%%%%%%
%%%%%%%%%%%%%%%%%%%%%%%%%%%%%%%%%%%%%%
%%%%%%%%%%%%%%%
 One can find details in section 1, \cite{IRO6}.

\subsection{4D case: variable rank}
\label{sect-2-5}
The most natural model operator corresponding to the canonical form
(\ref{7}) is $H^0+ H''$ with $H^0$ as above and $2H''=\xi_3^2+(\xi_4- x_3)^2$.
Then $H''$ is a movement integral. Therefore the dynamics is split into
dynamics in $(x',\xi')=(x_1,x_2,\xi_1,\xi_2)$ described above with potential $W=V-2E $ and the standard cyclotron movement with energy $E$ in $(x'',\xi'')=(x_3,x_4,\xi_3,\xi_4)$.

Situation actually is way more complicated: considering
$H^0 + (1+\alpha x_1)E$ we arrive to the 1-D potential
$\cV-1(1+\alpha x_1)$ and playing with $E$ and $\alpha$ one can kill the drift even for $k\gg 1$ leading to many periodic trajectories.

Consider canonical form (\ref{8})
which in polar coordinates in $(x_3,x_4)$ becomes
\begin{equation}
\sigma =
d \Bigl( (x_1-{\frac 1 2}\rho^2) dx_2 + (x_1-{\frac 1 4}\rho^2)\rho^2 d\theta)\Bigr).
\label{22}
\end{equation}

The most natural classical Hamiltonian corresponding to this form is
\begin{equation}
2H=\xi_1^2 + \bigl(\xi_2-\mu (x_1-{\frac 1 2}\rho^2)\bigr)^2+ \varrho^2+r^{-2}\bigl(\vartheta-\mu (x_1-{\frac 1 4}\rho^2)\rho^2\bigr)^2-1
\label{23}
\end{equation}
with $\varrho,\vartheta$ dual to $\rho,\theta$.

Note that $\xi_2$ and $\vartheta$ are movement integrals and therefore $x_1-{\frac 1 2}\rho^2$ is preserved modulo $O(\mu^{-1})$. Based on this one can prove that

\begin{itemize}
\item There is a cyclotronic movement with the angular velocity $\asymp \mu^{-1}$ in the normal direction to parabolloid $\{-x_1+{\frac 1 2}\rho^2={\frac 1 2}{\bar \rho}^2\}$
\item combined in the zone $\{|x_1|\le c\rho^2\}$ with the movement similar to one described in 2D case in $(\rho,\theta)$-coordinates (with $\{x_1=0\}$ now equivalent to $\{\rho={\bar \rho}\}$) on the surface of this ellipsoid
\item and also combined some movement along $x_2$;
\item I did not consider zone $\{|x_1|\ge c\rho^2\}$ since it was not needed
for the spectral asymptotics.
\end{itemize}
 One can find details in section 1, \cite{IRO8}.

\section{Quantum Dynamics}
\label{sect-3}
Microlocal canonical form (\emph{Birghoff normal form}) play a crucial role in the analysis of the quantum dynamics and spectral asymptotics.
\subsection{Canonical forms. I}
\label{sect-3-1}
In the case $d=2$ and a full-rank magnetic field  canonical form  of Magnetic Schr\"oding operator is (${\frac 1 2}$ of)
\begin{multline}
\omega_1 (x_1,\mu^{-1}hD_1) (h^2D_2^2+\mu^2 x_2^2) - W (x_1,\mu^{-1}hD_1)+\\
\sum_{m+k+l\ge2} a_{mkl} (x_1,\mu^{-1}hD_1) (h^2D_2^2+\mu^2 x_2^2)^m \mu^{2-2m-2k-l} h^l
\label{24}
\end{multline}
with $\omega_j = f_j\circ \Psi$, $W=V\circ \Psi$ with some map $\Psi$.
The first line is \emph{main part} of the canonical form.

In the case $d=3$ and a maximal-rank magnetic field microlocal canonical form  of Magnetic Schr\"oding operator is (${\frac 1 2}$) of
\begin{multline}
\omega_1 (x_1,x_2,\mu^{-1}hD_2) (h^2D_3^2+\mu^2 x_3^2) +h^2D_1^2 - W (x_1,x_2,\mu^{-1}hD_2)+\\
\sum_{m+n+k+l\ge2} a_{mnkl} (x_2,\mu^{-1}hD_2) (h^2D_3^2+\mu^2 x_3^2)^m D_1^n \times\\ \mu^{2-2m-2k-l-n} h^{l+n}
\label{25}
\end{multline}
Again, the first line is \emph{main part} of the canonical form.

In the case $d\ge 4 $ and a constant rank magnetic field microlocal canonical form  of Magnetic Schr\"oding operator
is of the similar type provided we can avoid some obstacles:

If $f_j$ have constant multiplicities (say, $f_j$  are simple for simplicity) then the main part is
\begin{multline}
\sum_{1\le j\le r} \omega_j(x',x'',\mu^{-1}hD'')(h^2D_{r+q+j}^2+\mu^2x_{r+q+j}^2) + h^2D'{}^2-\\
W (x',x'',\mu^{-1}hD'');
\label{26}
\end{multline}
where $x'=(x_1,\dots,x_q)$, $x''=(x_{q+1},\dots,x_{q+r})$, $2r=\rank F$, $q=d-2r$.

Next terms appear if one can avoid higher order \emph{resonances}: 
$\sum_j p_j f_j(x)=0$
with $p_j\in {\mathbb Z}$; $3\le \sum_j|p_j|$ is calleed \emph{the order of the resonance}.

After operator is reduced to the canonical form one can decompose functions as
\begin{equation}
u(x)=\sum_{\alpha \in {\mathbb Z}^{+r}} u_\alpha (x',x'')
\Upsilon_{p_1}(x_{r+q+1})\cdots \Upsilon_{p_r}(x_{d})
\label{27}
\end{equation}
where $\Upsilon $ are eigenfunctions of Harmonic oscillator $h^2D^2+\mu^2x^2$
(i.e. scaled Hermite functions).

Then as $2r=d$ one gets a family of $r$-dimensional $\mu^{-1}h$-PDOs and for $2r<d$ one  gets a family of $q$-dimensional Schr\"odinger operators with potentials which are $r$-dimensional $\mu^{-1}h$-PDOs.

The similar approach also works for $2$ and $4$-dimensional Schr\"odinger operators with the degenerate magnetic field of the types I considered before but only in the \emph{far outer zone} $\{\gamma(x)\Def|x_1|\gg \mu^{-1/\nu}\}$ and to this form operator is reduced in balls $B({\bar x}, {\frac 1 2}\gamma ({\bar x}))$.

\subsection{Canonical forms. II}
\label{sect-3-2}
As $d=2,4$ and magnetic field degenerates there is also a more global canonical form.

As $d=2$ in zone $\{|x_1|\ll 1\}$ this form is (after multiplication by some non-vanishing function)
\begin{equation}
h^2D_1^2 + (hD_2-\mu x_1^\nu /\nu)^2 -W(x) + \text{perturbation}
\label{28}
\end{equation}
with $W=V\phi^{-2/(\nu+1)}$ (if oiriginally $f_1\sim \phi \dist(x,\Sigma)^{\nu -1}$, $\Sigma=\{x_1=0\}$).

For $d=4$ one can separate a cyclotron part corresponding to the non-vanishing eigenvalue $f_2$; after this one gets a 3-dimensional second-order DO ($+$ perturbation) with the principal part which is the quadratic form of rank 2 and
a free term $V- (2\alpha +1)\mu h f_2$ where $\alpha\in {\mathbb Z}^+$ is a corresponding \emph{magnetic quantum number}.

This operator could be reduced to the form similar to (\ref{28})
(at least away from $\Lambda=\{x_3=x_4=0\}$); here
$W= \bigl(V- (2\alpha +1)\mu h f_2\bigr)\phi^{-2/3}$.

Close to $\Lambda$ but as $|x_1|\le C\rho^2$ one can get a similar form but with $\theta$ instead of $x_2$ and $\mu \rho$ instead of $\mu$.

\subsection{Periodic orbits}
\label{sect-3-3}
One can prove that semiclassical quantum dynamics follows the classical one long enough to recover sharp remainder estimates but the notion of \emph{periodic orbit} should be adjusted to reflect \emph{logarithmic uncertainty principle}
\begin{equation}
|\osc (x)|\cdot |\osc (D)|\ge C\hbar |\log \hbar|
\label{29}
\end{equation}
where $\hbar$ is effective Plank constant (it could be $h$ or $\mu^{-1}h$ or one of them scaled depending on the particular situation).

I need  a logarithm because I am interesting in the \emph{size of the box} outside of which function is negligible rather than in the mean quadratic deviation. Function $\exp ({-|x|^2/2\hbar })$ scaled shows shows that \emph{boxing} requires a logarithmic factor.

So instead of individual trajectories I consider their beams satisfying logarithmic uncertainty principle. One can see that the classical trajectory is not periodic but cannot say this about the semiclassical beam until much larger time, after few ``periods''.
\begin{figure}[h!]
\centering
\subfigure[ ]{
\includegraphics[width=0.4\textwidth]{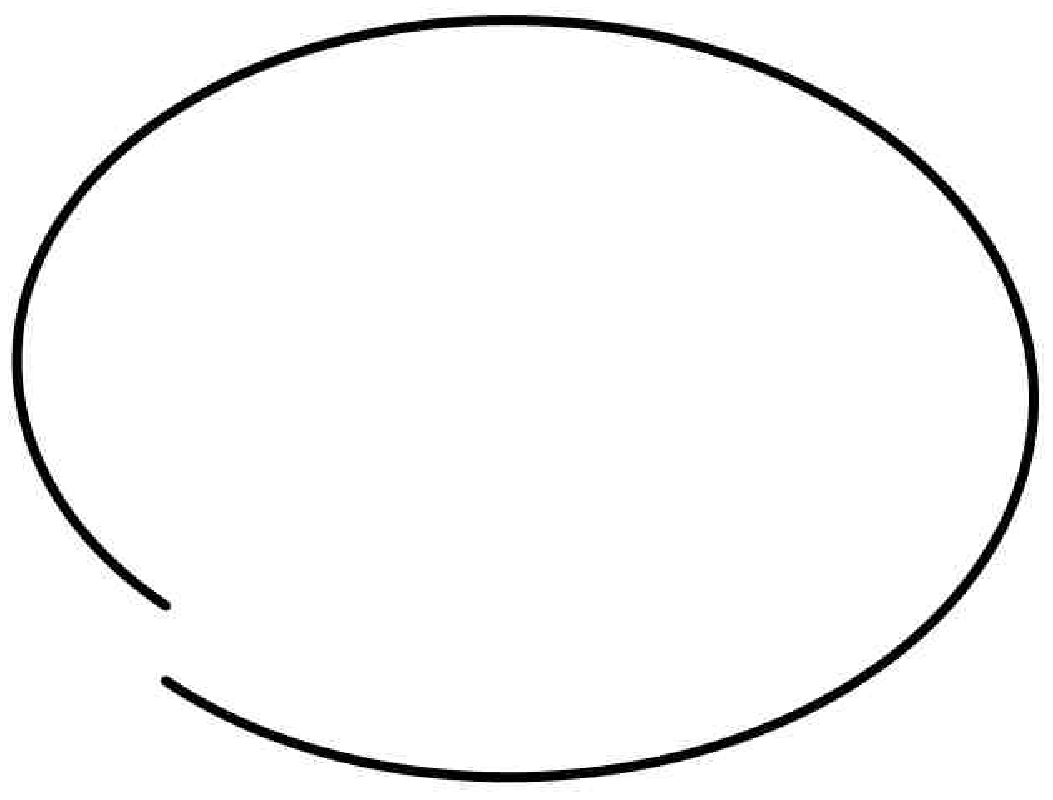}}
\subfigure[ ]{
\includegraphics[width=0.4\textwidth]{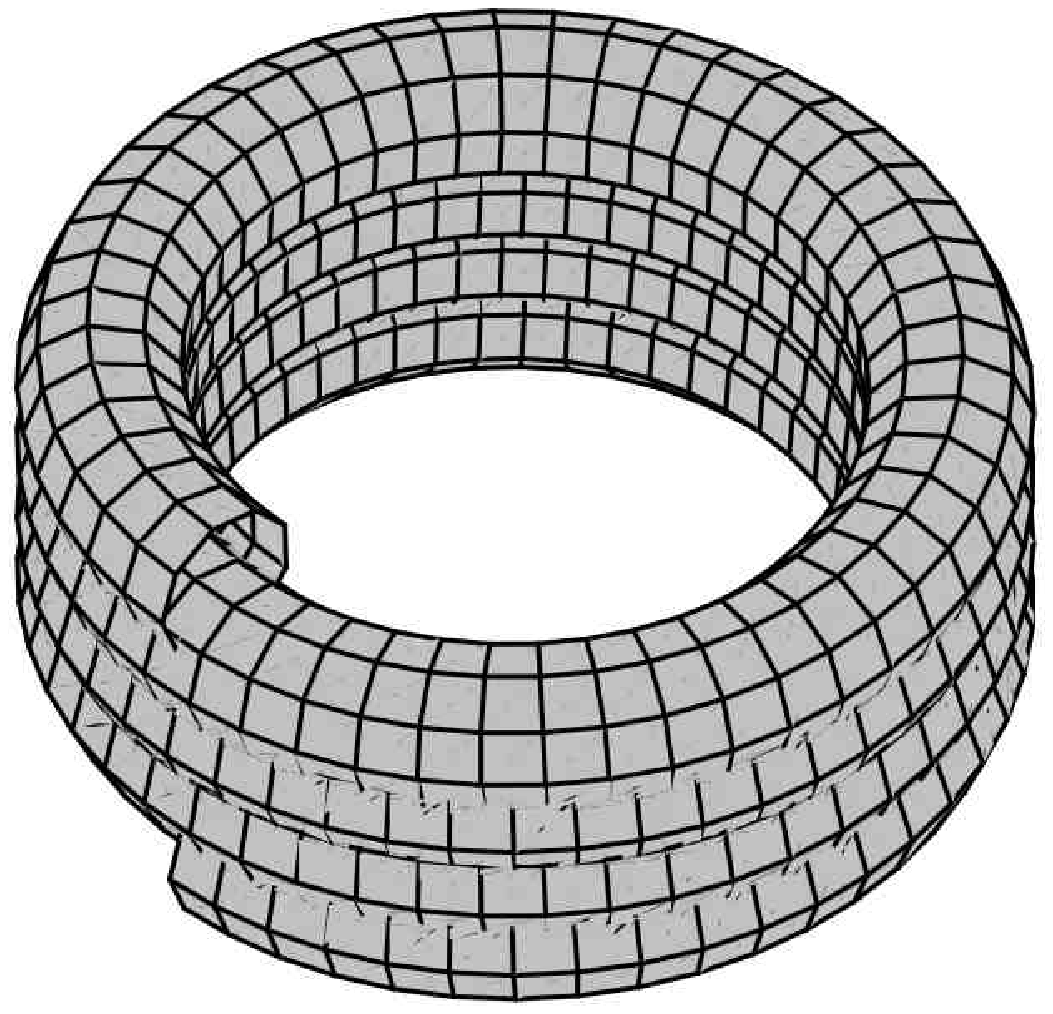}}
\caption{Classical (a) and semiclassical (b) periodicity.}
\end{figure}

\section{Spectral asymptotics}
\label{sect-4}
\subsection{Tauberian method}
\label{sect-4-1}
 \paragraph{Problem}
\label{par-4-1-1}
I am looking at asymptotics as $h\to +0$, $\mu \to +\infty$ of
\begin{equation}
\Gamma (Qe)(0)=\int \bigl(e(.,.,0)Q^t_y\bigr)_{x=y}dy = \Tr (QE(0))
\label{30}
\end{equation}
where $e(x,y,\tau)$ is the Schwartz kernel of the spectral projector
$E(\tau)$ of operator $H$ and $Q$ is a pseudo-differential operator, $Q^t$ means a dual operator.

As $Q=I$ we get $\Tr E(0)$ which is the number of negative eigenvalues of $H$
(and $+\infty$ if there is an essential spectrum of $H$ below $0$).

I hope to construct this expession (\ref{30}) with $Q=I$ from itself for elements of the partition of unity with self-adjoint $Q \ge 0$.

 \paragraph{Tauberian method}
\label{par-4-1-2}
\emph{Tauberian method Fourier} says that the \emph{main part} of $\Gamma (eQ_y^t)$ is given by expression
\begin{equation}
h^{-1}\int_{-\infty}^0 \Bigl(F_{t\to h^{-1}\tau}{\bar\chi}_T(t)\Gamma (uQ_y^t)\Bigr)\,d\tau
\label{31}
\end{equation}
while the \emph{remainder} does not exceed $C{\frac M T}+C'h^s$ where
\begin{equation}
M=M_T=\sup_{|\tau|\le \epsilon} |\Bigl(F_{t\to h^{-1}\tau}{\bar\chi}_T(t)\Gamma (uQ_y^t)\Bigr)|
\label{32}
\end{equation}
and $s$ is large, $C$ does not depend on $\epsilon, T,h,\mu$ and $s$ while $C'$ depends on $\epsilon>0$, $T$, $s$.

Here and below $u(x,y,t)$ is the Schwartz kernel of the propagator $e^{ih^{-1}tH}$,
${\bar\chi}\in C_0^\infty([-1,1])$ equal 1 at $[-{\frac 1 2}, {\frac 1 2}]$,
$\chi\in C_0^\infty([-1,1])$ equal 0 at $[-{\frac 1 2}, {\frac 1 2}]$,
$\chi_T(t)=\chi(t/T)$, $T>0$ and $F_{t\to h^{-1}\tau}$ is $h$-Fourier transform.

Actually this remainder estimate persists if one replaces $T$ by any larger number $T'$ \emph{only} in expression (\ref{31}).

So, I want to increase $T$ without (significantly) increasing $M_T$ in (\ref{32}).

\paragraph{Evil of periodic trajectories}
\label{par-4-1-3}
\emph{Microlocal analysis} says that if there are no periodic trajectories with periods in $[{\frac T 2},T]$ on energy levels in $[-2\epsilon,2\epsilon]$ then
\begin{equation}
\sup_{|\tau|\le \epsilon} |\Bigl(F_{t\to h^{-1}\tau}\chi_T(t)\Gamma (uQ_y^t)\Bigr)|\le C'h^s.
\label{33}
\end{equation}

Therefore if there are no periodic trajectories with periods in $[T ,T']$ on energy levels in $[-2\epsilon,2\epsilon]$, then one can retain
$T$ in (\ref{31}), $M_T$ in (\ref{32}) but the remainder estimate improves to
$C{\frac M {\emph{T'}}}+C'h^s$.

So, periodic trajectories are one of the main obstacles in getting a good remainder estimate. For example, if all trajectories are periodic with the period $T=T_\Pi$ then it can happen that $M_T \asymp TT_\Pi^{-1}M_{T_\Pi}$ as $T\ge T_\Pi$ and increasing $T$ does not bring any improvement.

For example, let $\mu\le 1$, $V\asymp 1$. Then there are no periodic trajectories with periods in $[T_0,T_1]$, $T_0=Ch|\log h|$ and $T_1=\epsilon$ because $\dist (x(t),x(0))\asymp T$ as $T\le T_1$ and this distance is observable as $T\ge T_0$.

Then $M\le C h^{-d}T_0= Ch^{1-d}|\log h|$ and the remainder estimate is\newline $O(h^{1-d}|\log h|)$. Actually one can get rid off this log factor in
 $M= Ch^{1-d}$ and the remainder estimate becomes $O(h^{1-d})$. This remainder estimate cannot be improved without geometric assumptions of the global nature.

Also, taking $T$ really small in (\ref{31}) allows us to calculate $u$ and (\ref{31}) by a crude \emph{successive approximation method} with an unperturbed operator $H$ having coefficients frozen as $x=y$.

As $\mu \ge 1$ the same arguments are true but $T_1= \epsilon \mu^{-1}$ and the remainder estimate is $O(\mu h^{1-d})$. This remainder estimate cannot be improved as $d=2$, $g^{jk}$, $f_1$ and $V$ are constant.

\emph{In our arguments\/}: because all trajectories are periodic (pure cyclotronic movement).

\emph{From direct calculations\/}: as domain is ${\mathbb R}^2$ all eigenvalues are
\emph{Landau levels} $(\alpha +{\frac 1 2})\mu h -{\frac 1 2} V$ of infinite multiplicity ($\alpha \in {\mathbb Z}^+$) and
\begin{equation}
e(x,x,\tau)= {\frac 1 {2\pi}} \sum_{n\ge 0}
\theta \bigl(2\tau +V-(2n+1)\mu h f\bigr) \mu h ^{-1}f\sqrt{g}
\label{34}
\end{equation}
with jumps $\asymp\mu h^{-1}$ at Landau levels.

However in many cases one can improve remainder estimate $O(\mu h^{1-d})$. The idea is to show that actually periodicity is broken.

From the point of view of applications one should take $Q$ with support (with respect to $x$) in ball \emph{$B(0,{\frac 1 2})$} (then rescaling arguments could be applied) and impose condition on operator \emph{only in the circle of light $B(0,1)$} with the self-adjointness being the \emph{only} condition outside of it.

So anything out of $B(0,1)$ is \emph{a dark territory} and we \emph{must} take
$T\le T^*$ which is the time for which trajectory which started from $\supp Q$
remains in $B(0,1)$. But we can chose the time direction and we can chose it for every beam individually.

Now, as $d=3$ the typical trajectory is non-periodic because of the free movement and
\begin{itemize}
\item one must take $T\le T^* \asymp 1$;
\item but for most of the trajectories one can take $T_1\asymp T^*$ retaining $T_0=Ch|\log h|$
\end{itemize}
and the remainder estimate is $O(h^{1-d})$ (under very mild assumptions).

\medskip
Let $d=2$ and $f_1$ do not vanish. Then \begin{itemize}
\item since the drift speed is $O(\mu^{-1})$
one can take $T \asymp \mu$;
\emph{under certain non-degeneracy conditions breaking periodicity of the cyclotronic movement} one can take $T_1=T^*$ retaining $T_0=Ch|\log h|$
\end{itemize}
and the remainder estimate is $O(\mu^{-1}h^{1-d})$.

When tamed, our worst enemy (periodic trajectories) could become our best friend!

\subsection{Results: Constant-rank case}
\label{sect-4-2}
 \paragraph{Results: ``Constant'' case}
 \label{par-4-2-1}
\begin{theorem}
\label{thm-1}
 Let $g^{jk},F_{jk}$ and $V$ be constant and domain be ${\mathbb R}^d$. Then
\begin{multline}
\cE_d^\MW (x,E)\Def \Omega_{d-2r}(2\pi )^{-d+r} \mu ^rh^{-d+r}\times\\
\sum _{\alpha \in {\mathbb Z}^{+r}}
\Bigl(2E +V - \sum_j (2\alpha_j +1) f_j\mu h -V\Bigr)_+^{{\frac d 2}-r}
f_1\cdots f_r\sqrt g
\label{35}
\end{multline}
where $\Omega_k$ is a volume of unit ball in ${\mathbb R}^k$.

\smallskip
\noindent
In particular, the spectrum is pure point iff $r=0$.
\end{theorem}
One can prove this theorem easily by direct calculations. \emph{Magnetic Weyl Expression\/} $cE_d^\MW (x,E)$ becomes our  candidate for the main part of asymptotics in the general case. 
 
\paragraph{Results: $d=2$}
 \label{par-4-2-2}
As $d=2$ formula (\ref{35}) provides a good approximation and the non-degeneracy condition below breaks periodicity and provides a good remainder estimate \cite{IRO16}:

\begin{theorem}\label{thm-2}  Let $d=2$ and $g^{jk}$, $F_{jk}$, $V$ be smooth in $B(0,1)$, $f_1$ non-vanishing there  and $\psi \in C_0^\infty (B(0,1))$. Let assume that
all critical values of $V/f_1$ are non-degenerate.
Then
\begin{equation}
|\int \bigl(e(x,x,0)- \cE_d^\MW (x,0)\bigr)\psi(x)\,dx|\le C\mu^{-1} h^{1-d}
\label{36}
\end{equation}
as $\mu \le ch^{-1}$.
\end{theorem}

\begin{remark}\label{rem-3} (i) If in the general case magnetic field spoils remainder estimate $O(\mu h^{1-d})$, but in the covered  case magnetic field improves it;

\smallskip
\noindent
(ii) Estimate (\ref{36}) holds in multidimensional full-rank case as well
but non-degeneracy condition is pretty complicated and is not generic \cite{IRO4};

\smallskip
\noindent
(iii) As $\epsilon_1\le \mu h\le c$ nondegeneracy condition changes; as $d=2$ it reads: $(2\alpha+1)\mu h $ is not a degenerate critical value of $V/f$ for
any $\alpha\in {\mathbb Z}^+$;

\smallskip
\noindent
(iv) As $f_1+\dots +f_r\ge \epsilon >0$, $e(x,x,0)$ is negligible and $\cE^\MW_d=0$ for $\mu \ge ch^{-1}$.
\end{remark}

\paragraph{Results: Non full-rank case}
\label{par-4-2-3}
Let us assume that $\rank F$ is constant but less than $d$. In this case remainder estimate cannot be better than $O(h^{1-d})$ but it also cannot be much worse \cite{IRO5}:

\begin{theorem}\label{thm-4}  Let $g^{jk}$, $F_{jk}$, $V$ be smooth in $B(0,1)$, $\rank (F)=2r$, $0<2r<d$ so $f_1$,\dots,$f_r$ do not vanish there  and $\psi \in C_0^\infty (B(0,1))$. Then

\smallskip
\noindent
{\rm (i)} As either $2r=d-1$ and some very mild non-degeneracy condition is fulfilled  or  $2r=d-2$  or  $\mu \le h^{\delta-1}$ with $\delta>0$ 
\begin{equation}
|\int \bigl(e(x,x,0)- \cE_d^\MW (x,0)\bigr)\psi(x)\,dx|\le C h^{1-d};
\label{37}
\end{equation}

\smallskip
\noindent
{\rm (ii)} As $2r=d-1$ the left-hand expression does not exceed $C\mu h^{2-\delta-d}+Ch^{1-d}$ with arbitrarily small $\delta>0$.
\end{theorem}

\subsection{Results: Degenerating 2D case}
\label{sect-4-3}
Consider case $d=2$, $f_1\asymp |x_1|^{\nu-1}$ with $\nu\ge2$ assuming
that
\begin{equation}
V\ge \epsilon_0>0.
\label{38}
\end{equation}
We consider $\epsilon$-vicinity of $\{x_1=0\}$ with small enough constant $\epsilon>0$.

Then in the \emph{outer zone} $\{{\bar\gamma}=C\mu^{-1/\nu}\le |x_1|\le \epsilon\}$ there is a drift with the speed $\mu^{-1}\gamma^{-\nu}$, this drift breaks periodicity and therefore contribution of the strip $\{|x_1|\asymp \gamma\}$ with $\gamma\in ({\bar\gamma},\epsilon)$ to the remainder estimate does not exceed $Ch^{1-d}\times\gamma \times \mu^{-1}\gamma^{-\nu}$ where the second factor is the width of the strip and the third one is the inverse ``control time''. Then the total contribution of the outer zone to the remainder estimate does not exceed the same expression as $\gamma={\bar\gamma}$ which is $C\mu^{-1/\nu}h^{1-d}$; this is our best shot.

In the \emph{inner zone} $\{|x_1|\le {\bar\gamma}\}$ or equivalently $\{|\xi_2|\le C_0\}$ the similar arguments work as long as $\rho\asymp |\xi_2-k^*V^{1/2}|\ge \epsilon$.

Otherwise there seems to be no drift to save the day.
However it is not that bad. Really, period in $x_1$ is $\asymp {\bar\gamma}$ and if
\begin{equation}
|\xi_2-k^*V^{1/2}|\asymp \rho
\label{39}
\end{equation}
the speed of the drift is $\asymp\rho$, the shift with respect to $x_2$ is 
$\asymp \rho {\bar\gamma}$ and in order to be observable it must satisfy logarithmic uncertainty principle $\rho{\bar\gamma}\times \rho \ge Ch|\log h|$ because characteristic scale in $\xi_2$ is $\rho$ now.
So, periodicity is broken provided
\begin{equation}
\rho\ge {\bar\rho}_1= C\bigl({\bar\gamma}^{-1}h|\log h|\bigr)^{\frac 1 2},
\label{40}
\end{equation}
which leaves us with much smaller \emph{periodic zone}
\begin{equation}
\cZ_{\rm per}=\bigl\{|\xi_2-k^*V^{1/2}|\le {\bar\rho}_1\bigr\}.
\label{41}
\end{equation}
And in this periodic zone picking up $T_1\asymp {\bar\gamma}$ we can derive remainder estimate $O(h^{-1}{\bar\rho}_1)$ which does not exceed our dream estimate $C\mu^{-1/\nu}h^{-1}$ provided ${\bar\rho}_1\le {\bar\gamma}$ or
\begin{equation}
\mu \le C(h|\log h|)^{-{\frac \nu 3}}.
\label{42}
\end{equation}

Actually for the general operator rather than the model one we need to assume that $\rho\ge C{\bar\gamma}$ but this does not spoil our dream estimate.

So we need to consider periodic zone defined by (\ref{41}) assuming that (\ref{42}) does not hold.

 \paragraph{Inside of Periodic Zone}
\label{par-4-3-1}
Even in the periodic zone $\cZ_{\rm per}$ periodicity of trajectories can be broken as $W=V\phi^{-2/(\nu+1)}|_{x_1=0}$ is ``variable enough'' which leads us \cite{IRO6,IRO7} to 

\begin{theorem}\label{thm-5}
Let $d=2$, $f_1=\phi (x) |x_1|^{\nu-1}$ with $\nu\ge2$ and condition $(\ref{38})$ be fulfilled. Then as $\psi \in C_0^\infty(B(0,1)\cap\{|x_1|\le \epsilon\}$
\begin{equation}
|\int \bigl(e(x,x,0)- \cE_d^\MW (x,0)\bigr)\psi(x)\,dx|\le  C(\mu^{-1/\nu}+\hbar^{(q+1)/2})h^{1-d}
\label{43}
\end{equation}
where here and below $\hbar= \mu^{1/\nu}h$, $q=0$ in the general case and $q=1$ under assumption ``$W$ does not have degenerate critical points''.
\end{theorem}

To improve this remainder estimate one should take in account the short periodic trajectories. Actually, periodicity of the trajectories close to them is broken but only after time $T_0=C\rho^{-2}h|\log h|$ (see our discussion in subsection \ref{sect-3-3}). Skipping details \cite{IRO6,IRO7} 

\begin{theorem}\label{thm-6}
Let $d=2$, $f_1=\phi (x) |x_1|^{\nu-1}$ with $\nu\ge2$ and condition $(\ref{38})$ be fulfilled. Then as $\psi \in C_0^\infty(B(0,1)\cap\{|x_1|\le \epsilon\}$
\begin{multline}
|\int \bigl(e(x,x,0)- \cE_d^\MW (x,0)\bigr)\psi(x)\,dx-\int \cE^\MW_{\rm corr}(x_2,0)\psi(0,x_2)\,dx_2|\le \\ C\mu^{-1/\nu}h^{1-d}
\label{44}
\end{multline}
provided either some very mild nondegeneracy condition is fulfilled  or  $\mu \le ch^{\delta-\nu}$ where $\cE^\MW_{\rm corr}(x_2,0)$ is defined below.
\end{theorem}

\begin{remark}\label{rem-7}
(i) In the case we are considering right now (and no other case considered in this article)
condition $f_1+\dots+f_r\ge \epsilon_0$ fails and therefore $e(x,x,0)$ is not negligible as $\mu \ge ch^{-1}$;

\smallskip
\noindent
(ii) On the other hand,  $e(x,x,0)$ is negligible as $\mu \ge ch^{-\nu}$;

\smallskip
\noindent
(iii) As $ch^{-1}\le \mu \le ch^{-\nu}$, $\cE^\MW(x,0)$ is supported in the strip
$\{|x_1|\le {\bar\gamma}_1\Def C_0(\mu h)^{-1/(\nu-1)}\}$ where
${\bar\gamma}_1\ge {\bar\gamma}=c\mu^{-1/\nu}$; therefore the main part of the spectral asymptotics (after integration) is of magnitude $(\mu h)^{-1/(\nu-1)}h^{-d}$.
\end{remark}

In  theorem \ref{thm-6} the correction term is defined by
\begin{equation}
\cE^\MW_{\rm corr}(x_2,\tau) = (2\pi h)^{-1}\int
{\bf n}_0(\tau; x_2, \xi_2,\hbar) \,d\xi_2 - \int \cE^\MW_0 (\tau; x_1,x_2,\hbar)\,dx_1
\label{45}
\end{equation}
where ${\bf n}_0$ is an eigenvalue counting function for an auxillary 1D-operator
\begin{equation}
{\bf a}_0(x_2,\xi_2,\hbar)={\frac 1 2}\Bigl(\hbar^2 D_1^2+
\bigl(\xi_2 - x_1^\nu/\nu\bigr)^2-W(x_2)\Bigr)
\label{46}
\end{equation}
and $\cE_0^\MW$ is Magnetic Weyl approximation for the related 2-dimensional operator.

Using Bohr-Sommerfeld approximation one can calculate eigenvalues
of ${\bf a}_0(x,\xi_2,\hbar)$ with $O(\hbar^s)$ precision and
$\cE^\MW_{\rm corr}(x_2,\tau)$ with $O(h^{-1}\hbar^s)$ precision.
In particular, modulo $O(h^{-1}\hbar)=O({\bar\gamma}^{-1})$
\begin{equation}
\cE^\MW_{\rm corr}(x_2,0) \equiv \varkappa h^{-1}\hbar^{\frac 1 2}
W^{{\frac 1 4}-{\frac 1 {4\nu}}} G\Bigl({\frac {S_0W^{{\frac 1 2}+{\frac 1 {2\nu}} } }{2\pi\hbar}}\Bigr)
\label{47}
\end{equation}
with some constants $\varkappa$ and $S_0$ and function $G$ defined by
\begin{equation}
G(t) = \int_\bR \Bigl(t+{\frac 1 2}\eta^2 - \bigl\lfloor t+{\frac 1 2}\eta^2 +{\frac 1 2}\bigr\rfloor \Bigr)\,d\eta
\label{48}
\end{equation}
with the converging integral in the right-hand expression.
One can prove easily that
\begin{equation}
G\not\equiv 0,\; G(t+1)=G(t),\;\int_0^1G(t)\,dt=0,\;
G \in C^{\frac 1 2}.
\label{49}
\end{equation}
This is one of examples of the short periodic trajectories really contributing to the asymptotics.

\subsection{Results: Degenerating 4D case}
\label{sect-4-4}
4D case is way more complicated than 2D one. But there are some good news: since $f_1+f_2\ge \epsilon$ in the generic case, we need to consider only $\mu \le ch^{-1}$.

The main difficulty in 4D case comes from the \emph{outer zone} $\{{\bar\gamma}=C\mu^{-1/\nu}\le |x_1|\le \epsilon\}$ because there could be short periodic trajectories. In other words: Landau level
$(2\alpha_1+1)\mu h f_1 + (2\alpha_2+1)\mu h f_2$ could be flat $0$.
Still in contrast to the very general case when this can happen for up to
$(\mu h)^{-1}$ pairs $\alpha \in {\mathbb Z}^{+2}$, in the assumptions of theorem below it can happen only for no more than $C$ pairs and I was able to prove \cite{IRO8}

\begin{theorem}\label{thm-8} Let $F$ is of Martinet-Roussarie type and
condition $(\ref{38})$ be fulfilled.
Then as $\psi$ is supported in $B(0,1)\cap\{|x_1|\le \epsilon\}$
\begin{equation}
|\int \bigl( e (x,x,0)-\cE_4^\MW (x,0)\bigr)\psi (x)\,dx| \le 
C\mu^{-1/2}h^{-3}+ C\mu^2 h^{-2}.
\label{50}
\end{equation}
\end{theorem}

One can improve this result under extra condition breaking flat Landau levels \cite{IRO8}:

\begin{theorem}\label{thm-9}
In frames of above theorem assume that \emph{$(V/f_2)_{x_1=0}$ does not have degenerate critical points}.
Then as $\psi$ is supported in $B(0,1)\cap\{|x_1|\le \epsilon\}$
\begin{equation}
|\int \bigl( e (x,x,0)-\cE_4^\MW (x,0)-\cE ^\MW_{{\rm corr}}(x,0)\bigr)\psi (x)\,dx| \le  C\mu^{-1/2}h^{-3}.
\label{51}
\end{equation}
Here $\cE^\MW_{\rm corr}=O(\mu^{5/4}h^{-3/2})$ is associated with \emph{periodic zone} $\{|x_1|\le c\mu^{-1/2}\}$, and is the sum of similar expressions in 2D case for $V_\beta=V-(2\beta+1)\mu hf_2$ with $\beta\in {\mathbb Z}^+$; locally all of them but one could be dropped.
\end{theorem}

\subsection{Results: Nondegenerating 4D case revisited}
\label{sect-4-5}
As I mentioned, even if magnetic field does not degenerate, non-degeneracy condition of \cite{IRO4} is not generic (it defines open but not everywhere dense set). Recently I was able to prove \cite{IRO9}

\begin{theorem}\label{thm-10} Let us consider 4D Schr\"odinger operator with non-degenerate generic magnetic field. Then

\smallskip
\noindent
{\rm (i)} For generic potential $V$ asymptotics holds woth the principal 
part as in $(\ref{51})$ and the remainder estimate $O(\mu^{-1}h^{-3})$;

\smallskip
\noindent
{\rm (ii)} For general potential $V$ asymptotics holds woth the principal 
part as in $(\ref{50})$ and the remainder estimate $O(\mu^{-1}h^{-3}+\mu^2h^{-2})$; correction term could be skipped with no penalty unless $h^{-1/3+\delta}\le \mu\le h^{-1/3-\delta}$ in which case it can be skipped with the penalty $O(h^{-8/3-\delta}$ where $\delta>0$ is arbitrarily small.
\end{theorem}

\hfill (added \today)

\smallskip
All results here are proven in the series of articles at

\centerline{\url{http://www.math.toronto.edu/ivrii/Research/Preprints.php}}

\input{IRO18.bbl}

\end{document}

%% file: IRO18.bbl
\bibliographystyle{alpha}

\providecommand{\bysame}{\leavevmode\hbox to3em{\hrulefill}\thinspace}

\vglue .06truein

%\hfill\hfill \emph{  \today \/}

\vglue .06truein

\begin{tabular}{rrl}
&{\hskip 220 pt} &Department of Mathematics,\cr
&&University of Toronto,\cr
&&40, St.George Str.,\cr
&&Toronto, Ontario M5S 2E4\cr
&&Canada\cr
&&ivrii@math.toronto.edu\cr
&&Fax: (416)978-4107\cr
\end{tabular}